\documentclass[final,10pt]{article}
\usepackage{amssymb,amsmath,amsfonts,showkeys,amsthm}

\everymath{\displaystyle}

\newtheorem{theorem}{Theorem}[section]
\newtheorem{lemma}[theorem]{Lemma}
\newtheorem{proposition}[theorem]{Proposition}
\newtheorem{corollary}[theorem]{Corollary}

\theoremstyle{definition}
\newtheorem{definition}[theorem]{Definition}

\theoremstyle{remark}
\newtheorem{remark}[theorem]{Remark}

\numberwithin{equation}{section}

\newcommand{\ds}{\displaystyle}
\newcommand{\hu}{\widehat{\textbf{u}}}
\newcommand{\hr}{\widehat{\rho}}
\newcommand{\rn}{\rho^n}
\newcommand{\un}{\textbf{u}^n}
\newcommand{\pn}{\mbox{\boldmath $\psi$}^n}
\newcommand{\real}{\mbox{I${\!}$R}}

\begin{document}

\title{Error bounds for semi-Galerkin approximations of nonhomogeneous
incompressible fluids}

\author{P. Braz e Silva\footnotemark[1]\ \footnotemark[2]
and M. A. Rojas-Medar\footnotemark[3] \footnotemark[4]}%
\date{}
\maketitle
\renewcommand{\thefootnote}{\fnsymbol{footnote}}
\footnotetext[1]{Departamento de Matem\'atica - Universidade Federal de Pernambuco \\
CEP 50740-540 Recife, PE, Brazil  ({\tt pablo@dmat.ufpe.br})}
\footnotetext[3]{Instituto de Matem\'atica, Estat\'istica e
Computa\c c\~ao Cient\'ifica - UNICAMP, Cx. Postal 6065, CEP
13083-970, Campinas, SP, Brazil ({\tt marko@ime.unicamp.br})}
\footnotetext[2]{P. Braz e Silva was supported for this work by FAPESP/Brazil, Grant \#02/13270-1
and is currently supported in part by CAPES/MECD-DGU Brazil/Spain, Grant \#2137-05-4.}
\footnotetext[4]{M. Rojas-Medar is partially supported by CNPq-Brazil, Grant \#301354/03-0,
CAPES/MECD-DGU Brazil/Spain, grant \#2137-05-4 and project BFM2003-06446-CO-01, Spain.}
\thanks{Math Subject Classification: 35Q30, 76M22, 65M15, 65M60}


\renewcommand{\thefootnote}{\arabic{footnote}}






\begin{abstract}
We consider spectral semi-Galerkin approximations
for the strong solutions of the nonhomogeneous Navier-Stokes equations. We
derive an optimal uniform in time error bound in the
$H^1$ norm for approximations of the velocity. We also derive an error estimate
for approximations of the density in some spaces $L^r$.
\end{abstract}

\maketitle

\section{Introduction}
Let $\Omega \subseteq \mathbb{R}^{n}, n =2$ or $3$, a
$C^{1,1}-$regular bounded domain, and $T > 0$. We are interested in the initial boundary
value problem
\begin{equation}
\left\{\begin{array}{ll} \label{eq1}
 \rho (\textbf{u}_t  +  \textbf{u} \cdot \nabla\textbf{u})-
\triangle \textbf{u}  + \nabla p = \rho {\bf f}  \; \mbox{ in } \; \Omega\times[0,T) ,
\\
  \mbox{div} \, \textbf{u}=0 \; \mbox{ in } \;  \Omega\times(0,T),
\\
 \rho_t + \textbf{u} \cdot \nabla \rho =  0 \; \mbox{ in } \; \Omega\times(0,T),\\
\textbf{u} = 0 \; \mbox{ on }\; \partial\Omega\times [0,T),  \\
\textbf{u}(x , 0)=\textbf{u}_0 (x), \; \forall \, x \in \Omega, \\
 \rho (x,0) = \rho_0(x), \; \forall \, x \in \Omega .
\end{array}\right.
\end{equation}
These are the equations of motion for nonhomogeneous
incompressible fluids.
The unknowns are
the velocity $\textbf{u}(x,t) \in \mathbb{R}^{n}$ of the
fluid, its density
$\rho(x,t)\in \mathbb{R}$ and the hydrostatic pressure $p(x,t)\in \mathbb{R}$.
The functions $\textbf{u}_0(x)$ and $\rho_0(x)$  are respectively the initial velocity
and initial density. The function $\textbf{f}(x,t)$ \ is the density by
unit of mass of the external force acting on the fluid. Here,
without loss of generality to our aim, the viscosity is considered to
be one.

Different methods have been used to study the existence, uniqueness and regularity of
solutions for problem (\ref{eq1})(see
\cite{antonzev}, \cite{boldrini1}, \cite{boldrini2}, \cite{kazhikov}
\cite{kim},  \cite{olga2}, \cite{JLLions2}, \cite{JLLions1}, \cite{PLLions}, \cite{okamoto},
\cite{salvi2}, \cite{simon1}). Here, we are interested in the spectral Semi-Galerkin method
applied to system (\ref{eq1}), the word spectral being used to indicate that the eigenfunctions
of the associated Stokes operator are used as a basis of approximations.
It is very important to derive error bounds for Galerin methods, due to the
wide application of these methods in numerical experiments. Even the case of
spectral Galerkin method may be used as a preparation and guide
for the more practical finite element Galerkin method.

A
systematic development of error bounds for the spectral
Galerkin method applied to the classical Navier-Stokes equations
was given in \cite{rautmann}. These error bounds are local in time in the
sense that they depend on functions which grow exponentially with
time. As observed in \cite{heywood1}, this is the best one may
expect without any assumptions about the stability of the solution
being approximated. Optimal uniform in time error
estimates for the velocity in the Dirichlet norm were also derived in \cite{heywood1},
assuming uniform boundedness in time of the $L^2-$norm of the
gradient of the velocity and exponential stability in the
Dirichlet norm of the solution.

For the variable density case, error
bounds were first obtained in \cite{salvi1}, where local in time error bounds were derived.
Moreover, a result of uniform in time error estimates in the
$L^2$ norm was also stated. This
last result, however, is not optimal. Moreover, it requires the
assumption $\textbf{u} \in L^\infty(0,T, \mathbf{H}^3(\Omega))$.
As pointed out in \cite{heywood2}, this assumption is pretty
restrictive, since it requires a global compatibility condition on
the initial data even for the classical Navier-Stokes equations. Error bounds were also
derived in \cite{boldrini4} without
explicitly assuming stability, but requiring exponential decay of
the external force field. This hypothesis though is very
restrictive as well, since gravitational forces do not satisfy it.

Here we derive error bounds assuming the
solution ($\textbf{u},\rho$) to be \textit{$p_0$-conditionally
asymptotically stable}, a notion to be defined in section \ref{description}.
The number $p_0$ is required to satisfy $6 \leq p_0 \leq\infty$, and is
related with the regularity of allowed perturbations of the density equation.
In \cite{heywood1}, a similar notion has been
used to treat the classical Navier-Stokes equations(see also
\cite{salvi1}). Here, we adapt it in the proper way to be used in
the variable density case. With this assumption, we obtain an
uniform in time optimal error bound in the Dirichlet norm for
the velocity. An error bound depending on time for the density
in some spaces $L^r$ is also derived.

In section \ref{preliminaries} we state some
preliminary useful results.
In section \ref{description} we describe the approximation method, the stability
notion to be used, and state
the main result. Section \ref{bounds} is dedicated to deriving
a priori estimates. Finally, we
present the proof of the main result in section \ref{proof}.

To simplify the notation, we denote by $C$
a generic finite positive constant depending only on $\Omega$ and
the other fixed parameters of the problem. As usual, it
may have different values in different expressions.
When necessary, we emphasize that the constants may have different values using
the notation $C_1$, $C_2$, and so on.

\section{Preliminaries} \label{preliminaries}

Throughout this work, we consider the usual Sobolev spaces
\[
W^{m,q}(D)=\{\,f \in L^q(D)\,, : ||\partial^\alpha f||_{L^q(D)} <
+ \infty\;,  |\alpha| \leq m\,\}\,,
\]
for a multi-index $\alpha $, a nonnegative integer $m$ and  $1
\leq q \leq +\infty $, where the domain is $D=\Omega$ or $D=\Omega \times
(0,T)$, $  0 < T \leq + \infty$, depending on the context. We write
$H^m(D):=W^{m,2}(D)$ and denote by $H^m_0(D)$ the closure of
$C^\infty_0(D)$ in $H^m(D)$. If $\;B\;$ is a Banach space, we
denote by $L^q([0,T]; B)$ the Banach space of $B$-valued functions
defined on the interval $[0,T]$ that are $L^q$-integrable in
Bochner's sense. Spaces of $\mathbb{R}^n$ valued functions, as well as their
elements, are represented
by bold face letters. We write
\[
\mathbf{C}^\infty_{0,\sigma}(\Omega):=\{\mathbf{v} \in
\mathbf{C}^\infty_0(\Omega) / \mbox{ div \bf{v} }=0 \quad
\mbox{in} \quad \Omega\},
\]
and denote by $\mathbf{H}$ and $\mathbf{V}$ the closure of
$\mathbf{C}^\infty_{0,\sigma}(\Omega)$ in $\mathbf{L}^2(\Omega)$ and
$\mathbf{H}^1(\Omega)$ respectively.

Throughout this work, the orthogonal projection from
$\mathbf{L}^2(\Omega)$ onto $\mathbf{H}$ is written as $P$. Thus,
the well known Stokes operator is  $-P \Delta$, defined on
$\mathbf{V} \cap \mathbf{H}^2(\Omega)$.
Its eigenfunctions and eigenvalues are denoted by
$\textbf{w}^k$ and $\lambda_k$ respectively. The usual
$L^2(\Omega)$ inner product and norm are respectively indicated by
$(\cdot ,\cdot )$, and $\|\cdot\|$.

It is well known that $\{\textbf{w}^k(x)\}^\infty_{k=1}$ form an
complete orthogonal system in the spaces $\mathbf{H}$,
$\mathbf{V}$ and $\mathbf{V} \cap \mathbf{H}^2(\Omega)$ equipped
with the usual inner products $(\textbf{u},\textbf{v}), (\nabla
\textbf{u},\nabla \textbf{v})$ and $(P\Delta \textbf{u},P\Delta
\textbf{v})$ respectively.

For each $k \in \mathbb{N}$, we denote by $P_k$ the orthogonal
projection from $\mathbf{L}^2(\Omega)$ onto $\mathbf{V}_k:=$
span$[\textbf{w}^1,...,\textbf{w}^k]$. For all
$\textbf{f},\textbf{\textbf{g}} \in \mathbf{L}^2(\Omega)$ and $k,m
\in \mathbb{N}$, it holds
\begin{itemize}
 \item[$(i)$]$(P_k \textbf{f},\textbf{g}) = (\textbf{f}, P_k \textbf{g})$,
\item[$(ii)$] $(P\textbf{f},\textbf{g}) = (\textbf{f}, P \textbf{g})$,
\item[$(iii)$] $((P_m-P_k) \textbf{f},\textbf{g}) = (\textbf{f},(P_m-P_k) \textbf{g})$,
\item[$(iv)$] $((P-P_k) \textbf{f},\textbf{g}) = (\textbf{f},(P-P_k) \textbf{g})$.
\end{itemize}
 The following bounds are useful to our ends.
\begin{lemma}[Rautmann\cite{rautmann}] \label{lema_rautmann}
If $\textbf{v} \in \mathbf{V}$, then
$$
\|\textbf{v}-P_k \textbf{v}\|^2 \leq \frac{1}{\lambda_{k+1}}
\|\nabla \textbf{v}\|^2.
$$
Moreover, if $\textbf{v} \in \mathbf{V} \cap
\mathbf{H}^2(\Omega)$, then
\[
\|\textbf{v}-P_k \textbf{v}\|^2 \leq \frac{1}{\lambda^2_{k+1}}
\|P\triangle \textbf{v}\|^2,
\]
\[
\|\nabla \textbf{v}-\nabla P_k \textbf{v}\|^2 \leq
\frac{1}{\lambda_{k+1}} \|P\triangle \textbf{v}\|^2.
\]
\end{lemma}
\begin{remark}
From Lemma \ref{lema_rautmann}, it follows that if $\textbf{f} \in
 \mathbf{H}^1(\Omega)$, then
$$
\|(I-P_k)P\textbf{f}\|^2 \leq \frac{1}{\lambda_{k+1}} \|\nabla P
\textbf{f}\|^2.
$$
Moreover, since $P : \mathbf{H}^1(\Omega) \rightarrow  \mathbf{
H}^1(\Omega)$ is a continuous operator\cite{wahl}, we have
\[
\|\nabla P\textbf{f}\|^2 \leq C\|\textbf{f}\|^2_{\mathbf{H}^1}.
\]
Thus, for all $\textbf{f} \in \mathbf{H}^1(\Omega)$, one has
\[
\|(I-P_k)P\textbf{f}\| \leq \frac{C}{\lambda_{k+1}}
\|\textbf{f}\|^2_{\mathbf{H}^1}
\]
Since $PP_k=P_kP=P_k$, one equivalently obtains
\[
\|P\textbf{f} - P_k\textbf{f}\|^2 \leq \frac{C}{\lambda_{k+1}}
\|\textbf{f}\|^2_{\mathbf{H}^1}.
\]
Moreover, the above relations also hold
if one replaces $P$ by any $P_m$, $m
>k$.
Analogously, one may check that
$$
\|(I-P_k)P\textbf{f}\|^2 \leq \frac{C}{\lambda^2_{k+1}}
\|\textbf{f}\|^2_{\mathbf{H}^2}
$$
for all $\textbf{f} \in \mathbf{H}^2(\Omega)$.
We also have, as an easy consequence
of the $L^2$-orthogonality of the functions
$\{\textbf{w}^k\}^\infty_{k=1}$, the following: Let $m>k$, $m,k
\in \mathbb{N}$, $\textbf{f} \in \mathbf{L}^2(\Omega)$ and
$\textbf{v}^m \in \mathbf{V}_m$, $ \textbf{v}^k \in \mathbf{V}_k$.
Then
\[
((P_m-P_k)\textbf{f}, \textbf{v}^m-\textbf{v}^k) =
(\textbf{f},(I-P_k)\textbf{v}^m).
\]
\end{remark}
\section{Stability concept and main result} \label{description}

We consider problem (\ref{eq1}) for all time $t\geq 0$, and
suppose the data given therein satisfy
\begin{eqnarray}
   & & \displaystyle \textbf{u}_0 \in \mathbf{V} \cap \mathbf{H}^2 (\Omega ) , \label{eq2} \\
    & & \sup_{t \geq 0} \| \textbf{f}\|_{\mathbf{H}^1} < \infty \; ; \;
    \sup_{t \geq 0} \|\textbf{f}_t \| < \infty ,  \label{eq3} \\
    & & \rho_0 \in C^1 (\bar{\Omega}) \; ; \; 0 < \alpha \leq \rho_0 \leq \beta , \label{eq3.1}
\end{eqnarray}
where $\alpha$ and $\beta$ are constants. We also suppose
that there exists $M > 0 $ such that the solution $( \textbf{u} ,
\rho )$ of (\ref{eq1}) satisfies
\begin{equation}
\| \nabla \textbf{u}(t) \| \leq M \; ; \; \forall \, t \geq 0 .
\label{eq4}
\end{equation}
If $n = 2$, then conditions (\ref{eq2}) and (\ref{eq3})
imply that (\ref{eq4}) holds. If $n = 3$, then inequality (\ref{eq4})
holds for $\textbf{f}$ and $\textbf{u}_0$ sufficiently small (see
\cite{boldrini2}). For a given $p_0$, $6\leq p_0
\leq\infty$, we also assume $(\textbf{u} , \rho )$ to be {\it
$p_0$-conditionally asymptotically stable} (see
\cite{heywood1,salvi1} for similar notions of stability). To define
this notion of stability, we first define perturbations of system
(\ref{eq1}). The pair $( \mbox{\boldmath $\mbox{\boldmath
$\xi$}$} (x , t) , \eta (x , t)) $, defined for $t\geq
t_0\geq 0$, is called  a perturbation of $( \textbf{u}, \rho)$ at time $t_0$ if $(\hu
:= \textbf{u} + \mbox{\boldmath $\mbox{\boldmath $\xi$}$}, \hr :=
\rho + \eta)$ is a solution of (\ref{eq1}), with $\mbox{\boldmath
$\mbox{\boldmath $\xi$}$}|_{\partial \Omega} = 0$. Therefore,
setting $\mbox{\boldmath $\mbox{\boldmath $\xi$}$}_0 :=
\mbox{\boldmath $\mbox{\boldmath $\xi$}$} (\cdot , t_0)$, $\eta_0 :=
\eta( \cdot , t_0)$, the pair $(\mbox{\boldmath $\mbox{\boldmath
$\xi$}$} , \hr)$ is a solution of the initial boundary value problem
\begin{equation}  \label{eq5}
    \left\{ \begin{array}{l}
        \displaystyle \hr \mbox{\boldmath $\mbox{\boldmath $\xi$}$}_t + \hr (\textbf{u}\cdot \nabla )\mbox{\boldmath $\mbox{\boldmath $\xi$}$} + \hr (\mbox{\boldmath $\mbox{\boldmath $\xi$}$} \cdot \nabla ) \textbf{u}
        + \hr (\mbox{\boldmath $\mbox{\boldmath $\xi$}$}\cdot \nabla )\mbox{\boldmath $\xi$}
         + \nabla q =  \\ \hspace{2.8cm} = \Delta \mbox{\boldmath $\xi$} + (\rho - \hr ) ( \textbf{u}_t +
         (\textbf{u} \cdot \nabla) \textbf{u} -\textbf{f} )\; \; \mbox{in} \;\;  \Omega \times (t_0 , \infty),\\
        \displaystyle \hr_t + ((\mbox{\boldmath $\xi$} + \textbf{u})\cdot \nabla )\hr = 0
        \;\; \mbox{in} \;\; \Omega \times (t_0 , \infty),\\
        \displaystyle \nabla \cdot \mbox{\boldmath $\xi$} = 0 \;\; \mbox{in} \;\; \Omega \times (t_0 , \infty) , \\
        \mbox{\boldmath $\xi$}|_{\partial \Omega} = 0, \\
        \mbox{\boldmath $\xi$}(x, t_0) = \mbox{\boldmath $\xi$}_0 (x), \; \forall \, x \in \Omega,  \\
        \hr (x , t_0) = \rho (x , t_0) + \eta_0 (x), \; \forall \, x \in \Omega .
     \end{array} \right.
\end{equation}
Now, for a given $p_0$, $6\leq p_0\leq \infty$, we define the
concept of $p_0$-conditional asymptotic stability.
\begin{definition} \label{def1}
  The pair $( \textbf{\textbf{u}} , \rho )$ is said to be $p_0$-conditionally asymptotically
  stable if for
  all $t_0\geq 0$ there exist positive numbers $A$, $B$, $\delta$, $M_1 $, $M_2$ and a
  continuous decreasing function $F: [0 , \infty ) \rightarrow \real^+ $,
  $F(0) = 1$, $\ds \lim_{t\rightarrow \infty} F(t) = 0$
   such that, for all $\mbox{\boldmath $\xi$}_0 \in V \cap H^2(\Omega)$,
   $\eta_0 \in L^\infty(\Omega) \cap W^{1, p_0}(\Omega)$,
   satisfying
   $\| \nabla \mbox{\boldmath $\xi$}_0 \| < \delta$, $\| P\Delta \mbox{\boldmath $\xi$}_0 \| < A$, $\|  \eta_0 \|_{L^\infty} < B$,
    problem (\ref{eq5}) is uniquely solvable with
   \begin{eqnarray*}
         & & \mbox{\boldmath $\xi$} \in L^2_{loc} ( [t_0 , \infty ) ; \mathbf{V} \cap \mathbf{H}^2 (\Omega) ), \\
         & & \mbox{\boldmath $\xi$}_t \in L^2_{loc} ( [t_0 , \infty ) ; \mathbf{H}^1 (\Omega) ), \\
         & & \eta \in L^\infty ( [t_0 , \infty ) ; L^\infty(\Omega) \cap W^{1, p_0}(\Omega)).
   \end{eqnarray*}
   Moreover,
   \begin{eqnarray}
        & & \| \nabla \eta (\cdot , t )\|_{L^{p_0}} \leq M_1 \:  , \, \forall \,  t \geq t_0, \label{eq6} \\
         & & \| \nabla \mbox{\boldmath $\xi$} (\cdot , t) \| \leq M_2 \| \nabla \mbox{\boldmath $\xi$}_0 \| F (t - t_0 ) \, ,  \: \forall
         \, t \geq t_0. \label{eq7}
   \end{eqnarray}
\end{definition}
\begin{remark}
We use a general function $F(t)$ in
Definition \ref{def1} just to stress out that the results here do
not require an exponential decay rate.
\end{remark}
The solution of problem (\ref{eq1}) can be obtained through a
spectral semi-Galerkin approximation, that is, a spectral Galerkin
approximation
$$
 \textbf{u}^n ( x , t ) = \sum_{k =1}^{n} C_{kn} (t) \textbf{w}^k (x)
$$
for the velocity $\textbf{u}$, uniquely determined by
\begin{eqnarray}
  & & ( \rn \un_t , \mbox{\boldmath $\phi$}^n) + (\rn \un \cdot \nabla \un , \mbox{\boldmath $\phi$}^n ) + ( \nabla \un , \mbox{\boldmath $\phi$}^n)
  = (\rn \textbf{f} , \mbox{\boldmath $\phi$}^n), \; \,  t\geq 0, \label{eq8} \\
  & & (\un ( x , 0) - \textbf{u}_0 (x) , \mbox{\boldmath $\phi$}^n ) = 0 \label{eq8.1},
\end{eqnarray}
 for all $\mbox{\boldmath $\phi$}^n$ of the form $\mbox{\boldmath $\phi$}^n ( x) = \sum_{k=1}^{n} d_k \textbf{w}^k (x)$,
 and an infinite dimensional approximation $\rho^n$ for the density, solution of
\begin{equation} \label{eq9}
\begin{array}{rcl}
 & &  \rn_t + \un \cdot \nabla \rn = 0, \\
 & & \rn (0) = \rho_0 .
\end{array}
\end{equation}
It can be proved that $(\textbf{u}^n , \rho^n)$ converges in an
appropriate sense to $(\textbf{u} , \rho )$, solution of
(\ref{eq1}). Our main result is
\begin{theorem}\label{theorem1}
  Suppose $( \textbf{\textbf{u}} , \rho )$ to be
  $p_0$-conditionally asymptotically stable, for some $p_0$, $6 \leq p_0 \leq \infty$.
  Then, there exist constants $N \in \mathbb{N}$
  and $C \geq 0$ such that if $n \geq N$ then, for all $t\geq 0$,
 \begin{equation}
    \| \nabla \textbf{u} (\cdot , t ) - \nabla \un (\cdot , t) \| \leq \frac{C}{(\lambda_{n+1})^{\frac{1}{2}}}. \label{eq10}
  \end{equation}
Moreover, if $6 \leq p_0 < \infty$, then
\begin{eqnarray}
& & \| \rho (\cdot , t ) -  \rn (\cdot , t) \|_{L^r}
    \leq \frac{C t}{(\lambda_{n+1})^{\frac{1}{2}}} \;\: , \:\;  2 \leq r \leq \frac{6p_0}{6+p_0}, \label{eq11}
\end{eqnarray}
and if $p_0 = \infty$, then
\begin{equation}\label{eq11,1}
 \| \rho (\cdot , t ) -  \rn (\cdot , t) \|_{L^r}
    \leq \frac{C t}{(\lambda_{n+1})^{\frac{1}{2}}} \;\: , \:\;  2 \leq r \leq 6.
\end{equation}
 The constants $N$, $C$, depend only on the domain,
  on the norms of the data in (\ref{eq2}), (\ref{eq3}) and on the constants introduced in
  (\ref{eq4}) and definition \ref{def1}.
\end{theorem}

\section{A priori estimates} \label{bounds}

We first state a general simple result that will be used later on. A proof is given
in the Appendix.
\begin{lemma}\label{lemma11}
 Let $h(t)$ be an integrable nonnegative function. Suppose
there exist nonnegative constants $a_1$, $a_2$ satisfying
 $$
 \displaystyle \int_{t_0}^t h (\tau) d \tau \leq a_1 (t - t_0) + a_2,
 $$
 for all $t$, $t_0$ with $0\leq t_0 \leq t$.
Then,
$$
  \displaystyle \sup_{t \geq 0} \left(  e^{-t}\int_{0}^{t} e^{\tau} h (\tau ) d \tau\right) < \infty .
$$
\end{lemma}
For $\textbf{u}$, solution of (\ref{eq1}), and the perturbations
$\mbox{\boldmath $\xi$}$, we have:
 \begin{lemma}\label{lemma1}
    Given $\epsilon >0$, the bounds
 \begin{eqnarray}
  & & \frac{1}{2} \frac{d}{dt} \| \nabla \textbf{u}(\cdot , t) \|^2 +
  \frac{\epsilon}{8} \| P \Delta \textbf{u}(\cdot , t)\|^2
  + \frac{1}{4} \| \rho^{\frac{1}{2}} \textbf{u}_t (\cdot , t)\|^2
  \leq \label{eq12} \\
  & & \hspace{7cm}C\left( \| \textbf{f}(\cdot , t) \|^2 + \| \nabla \textbf{u}(\cdot , t) \|^6\right), \nonumber \\
& & \| \mbox{\boldmath $\xi$}_t (\cdot , t)\|^2 +
 \frac{1}{2} \frac{d}{dt} \| \nabla \mbox{\boldmath $\xi$} (\cdot , t)\|^2
 + \frac{\epsilon}{42} \| P \Delta \mbox{\boldmath $\xi$} (\cdot , t)\|^2
  \leq \label{eq13} \\ & & C \left( \| \nabla \mbox{\boldmath $\xi$} (\cdot , t)\|^2
 \| P \Delta \textbf{u}(\cdot , t) \|^2
    \mbox{} + \| \nabla \mbox{\boldmath $\xi$}(\cdot , t) \|^{6}
  + \|\textbf{u}_t(\cdot , t)\|^2 + \| P\Delta \textbf{u}(\cdot , t)\|^2 \right),\nonumber
\end{eqnarray}
hold for all $t\geq 0 $.
 \end{lemma}
\begin{proof} Inequality (\ref{eq12}) was proved in \cite{simon1}. Inequality
(\ref{eq13}) can be proved in a completely analogous way.
\end{proof}
\begin{corollary}\label{cor1}
For all $t\geq t_0$, we have
  \begin{eqnarray}
   & & \displaystyle \int_{t_0}^t \| P \Delta \textbf{u}(\cdot , \tau) \|^2 d\tau \leq C + C(t - t_0) ,\label{eq14} \\
    & & \displaystyle \int_{t_0}^t \| \textbf{u}_t (\cdot , \tau)\|^2 d\tau \leq C + C(t - t_0) , \label{eq15} \\
    & & \displaystyle \int_{t_0}^t \| P \Delta \mbox{\boldmath $\xi$} (\cdot , \tau)\|^2d\tau
    \leq C + C(t - t_0) ,\label{eq16} \\
    & & \displaystyle \int_{t_0}^t \| \mbox{\boldmath $\xi$}_t (\cdot , \tau)\|^2d\tau \leq C + C(t - t_0) . \label{eq17}
  \end{eqnarray}
  Moreover, combining inequalities (\ref{eq14}), (\ref{eq15}), (\ref{eq16}) and (\ref{eq17})
  with lemma \ref{lemma11}, one gets
  \begin{eqnarray}
   & & \displaystyle \sup_{t \geq 0}
    e^{-t} \int_{0}^{t} e^{\tau} \| \textbf{u}_t (\cdot , \tau) \|^2
   d \tau < \infty , \label{eq18}\\
   & & \ds \sup_{t \geq 0} e^{-t} \int_{0}^{t} e^{\tau} \| P \Delta \textbf{u} (\cdot , \tau)\|^2 d \tau
   < \infty, \label{eq19} \\
   & & \displaystyle \sup_{t \geq 0}
    e^{-t} \int_{0}^{t} e^{\tau} \| \mbox{\boldmath $\xi$}_t (\cdot , \tau) \|^2
   d \tau < \infty , \label{eq18.1}\\
   & & \ds \sup_{t \geq 0} e^{-t} \int_{0}^{t} e^{\tau} \| P \Delta \mbox{\boldmath $\xi$} (\cdot , \tau)\|^2 d \tau
   < \infty \label{eq19.1}.
 \end{eqnarray}
\end{corollary}
The following lemma states some bounds for $\textbf{u}$ which are very
important for our later arguments.
\begin{lemma} \label{lemma4}
We have
  \begin{eqnarray}
    &  \ds \sup_{t\geq 0} \|  \textbf{u}_t (\cdot , t) \|
    < \infty , & \label{eq22}\\
     &  \ds \sup_{t\geq 0} \|  P\Delta \textbf{u}  (\cdot , t) \|
   <  \infty , &   \label{eq23} \\
     &  \ds \sup_{t\geq 0} e^{-t} \int_{0}^{t} e^{\tau} \|\nabla \textbf{u}_t (\cdot , \tau)\|^2 d \tau
   <  \infty  .&  \label{eq24}
  \end{eqnarray}
\end{lemma}
\begin{proof} We begin by proving (\ref{eq23}), supposing
(\ref{eq22}) holds.
 Setting $\textbf{v} = - P \Delta \textbf{u}$ in the weak formulation of problem
 (\ref{eq1}), we get
 $$
  - ( \rho \textbf{u}_t , P \Delta \textbf{u} )  - (\rho \textbf{u} \cdot \nabla \textbf{u} , P \Delta \textbf{u} )
  + \| P \Delta \textbf{u} \|^2 = - ( \rho \textbf{f} , P \Delta \textbf{u} ).
 $$
Thus,
\begin{eqnarray*}
  \| P \Delta \textbf{u} \| &  \leq & \|\rho \textbf{u}_t \| + \|\rho \textbf{u} \cdot \nabla \textbf{u} \|
  + \|\rho \textbf{f} \|   \\
  & \leq & \beta \| \textbf{u}_t \| + \beta \| \textbf{u} \|_{L^6} \|\nabla \textbf{u} \|_{L^3} + \beta \| \textbf{f} \| \\
  & \leq & \beta \| \textbf{u}_t \| +
  C\beta
  \| P \Delta \textbf{u}\|^{\frac{1}{2}}\| \nabla \textbf{u}\|^{\frac{3}{2}}
  + \beta \| \textbf{f} \| \\
  & \leq & \ds \beta \left( \| \textbf{u}_t \| + \| \textbf{f} \|\right)
  + \frac{(C\beta)^2}{2} \|\nabla \textbf{u} \|^3 +
  \frac{1}{2} \| P\Delta \textbf{u} \| \\
  & \leq & 2 \beta \| \textbf{u}_t \|
  + 2\beta \| \textbf{f}\| + (C\beta)^2 M^3 +
  \frac{1}{2} \| P\Delta \textbf{u} \|.
\end{eqnarray*}
Therefore, by (\ref{eq22}) and (\ref{eq3}), we have
\begin{eqnarray*}
  \ds \sup_{t \geq 0} \| P \Delta \textbf{u} (\cdot , t)\|&  \leq & \ds
  2 \beta \sup_{t \geq 0} \| \textbf{u}_t (\cdot , t)\|
  + 2\beta \sup_{t \geq 0}\| \textbf{f} (\cdot , t)\|
  + (C\beta)^2 M^3  < \infty ,
\end{eqnarray*}
which proves (\ref{eq23}). To prove (\ref{eq22}) and (\ref{eq24}),
differentiate the weak formulation of problem (\ref{eq1}), and set
$\textbf{v} = \textbf{u}_t$ to get
$$
 ( \rho_t \textbf{u}_t , \textbf{u}_t) + ( \rho \textbf{u}_{tt} , \textbf{u}_t )
 + ( (\rho (\textbf{u} \cdot \nabla) \textbf{u})_t , \textbf{u}_t)
 + (\nabla \textbf{u}_t , \nabla \textbf{u}_t)
 = ((\rho \textbf{f})_t , \textbf{u}_t) .
$$
Therefore,
\begin{eqnarray}
 \frac{1}{2}\frac{d}{dt} \|\sqrt{\rho} \textbf{u}_t  \|^2 + \| \nabla \textbf{u}_t \|^2
 & = &
 -\frac{1}{2} \int_{\Omega}  \rho \textbf{u} \cdot \nabla (\textbf{u}_t \cdot \textbf{u}_t)\, dx
  \nonumber \\
 & & \mbox{} - ( (\rho (\textbf{u} \cdot \nabla) \textbf{u})_t , \textbf{u}_t)
 + ((\rho \textbf{f})_t , \textbf{u}_t) . \label{eq25}
\end{eqnarray}
Now, bound each term on the right hand side of (\ref{eq25}) as follows:
\begin{eqnarray*}
\left| -\frac{1}{2} \int_{\Omega} \rho \textbf{u} \cdot \nabla
(\textbf{u}_t \cdot \textbf{u} _t)\, dx \right| &
\leq  & \| \rho \|_{\infty} \| \textbf{u} \|_{L^4} \| \textbf{u}_t\|_{L^4} \| \nabla \textbf{u}_t \| \\
 & \leq & C \beta \| \nabla \textbf{u} \| \| \nabla \textbf{u}_t \| \{ \| \textbf{u}_t \|^{\frac{1}{4}} \| \nabla \textbf{u}_t \|^{\frac{3}{4}} \} \\
 & = & C \beta \|\nabla \textbf{u} \| \| \textbf{u}_t \|^{\frac{1}{4}} \| \nabla \textbf{u}_t \|^{\frac{7}{4}} \\
 & \leq &
 C_{\epsilon} (C \beta \|\nabla \textbf{u} \|)^8 \| \textbf{u}_t\|^2
 + \epsilon \| \nabla \textbf{u}_t \|^2 ,
\end{eqnarray*}
\begin{eqnarray*}
 | (\rho \textbf{f}_t , \textbf{u}_t ) | & \leq & \frac{\beta^2}{2} \| \textbf{f}_t \|^2 + \frac{1}{2} \|\textbf{u}_t \|^2 ,
\end{eqnarray*}
\begin{eqnarray*}
 | (\rho_t \textbf{f} , \textbf{u}_t ) | & = & \left| \int_{\Omega} \rho_t \textbf{f}\cdot \textbf{u}_t  \, dx \right|=
 \left|- \int_{\Omega} \mbox{div}(\rho \textbf{u}) \textbf{f} \cdot \textbf{u}_t \, dx\right| =
 \left|- \int_{\Omega} \rho \textbf{u}\cdot \nabla( \textbf{f} \cdot \textbf{u}_t) \, dx\right| \\
 & \leq & \| \rho \|_{L^{\infty}} \| \textbf{u} \| \| \nabla (\textbf{f} \cdot \textbf{u}_t) \| \leq
 C_{\epsilon} ( C \beta \| \nabla \textbf{u}\| \| \nabla \textbf{f}\|)^2 + \epsilon \| \nabla \textbf{u}_t \|^2 \\
 & \leq & C_{\epsilon} ( C \beta M \| \nabla \textbf{f}\|)^2 + \epsilon \| \nabla \textbf{u}_t \|^2,
\end{eqnarray*}
\begin{eqnarray*}
 | (\rho (\textbf{u}_t\cdot \nabla)\textbf{u} , \textbf{u}_t) | & \leq & \| \rho \|_{L^{\infty}} \| \textbf{u}_t \|_{L^4}^2 \|\nabla \textbf{u}\|
 \leq C \beta \|\nabla \textbf{u}\| \| \textbf{u}_t\|^{\frac{1}{2}}\| \nabla \textbf{u}_t\|^{\frac{3}{2}} \\
 & \leq & C_{\epsilon} ( C \beta \| \nabla \textbf{u}\| \|\textbf{u}_t\|^{\frac{1}{2}})^4 + \epsilon \| \nabla \textbf{u}_t \|^2
 \leq C_{\epsilon} C^4 \beta^4 M^4 \|\textbf{u}_t\|^2 + \epsilon \| \nabla \textbf{u}_t \|^2,
\end{eqnarray*}
\begin{eqnarray*}
 | (\rho (\textbf{u}\cdot \nabla)\textbf{u}_t , \textbf{u}_t) | & \leq & \| \rho \|_{L^{\infty}}
 \| \textbf{u}_t \|^{\frac{1}{4}} \| \nabla \textbf{u}\|\| \nabla \textbf{u}_t\|^{\frac{7}{4}}
 \leq C_{\epsilon} ( C \beta M )^{8} \| \textbf{u}_t \|^2 + \epsilon \| \nabla \textbf{u}_t \|^2,
\end{eqnarray*}
\begin{eqnarray*}
 | (\rho_t (\textbf{u}\cdot \nabla)\textbf{u} , \textbf{u}_t) |   & = &
 | \int_{\Omega} \rho_t (\textbf{u}\cdot \nabla)\textbf{u} \cdot \textbf{u}_t |  =
 | - \int_{\Omega} \mbox{div}( \rho \textbf{u})  (\textbf{u}\cdot \nabla)\textbf{u} \cdot \textbf{u}_t | \\
 & \leq & C C_{\epsilon} \beta^2 M^4 \| P \Delta \textbf{u}\|^2 + 4 \epsilon \| \nabla \textbf{u}_t \|^2 .
\end{eqnarray*}
Therefore, from equation (\ref{eq25}) one obtains
\begin{equation} \label{eq26}
\begin{array} {l}
\frac{1}{2}\frac{d}{dt} \|\sqrt{\rho} \textbf{u}_t  \|^2 + \|
\nabla \textbf{u}_t \|^2 \leq
C_{\epsilon} (C \beta M)^8 \| \textbf{u}_t \|^2 +
\frac{\beta^2}{2} \| \textbf{f}_t \|^2 +
C_{\epsilon} ( C \beta M \| \nabla \textbf{f}\|)^2  \\
\mbox{} + C_{\epsilon} C^4 \beta^4 M^4 \|\textbf{u}_t\|^2 +
C
C_{\epsilon} \beta^2 M^4 \| P \Delta \textbf{u}\|^2 + \left(8
\epsilon + \frac{1}{2}\right) \| \nabla \textbf{u}_t \|^2 .
 \end{array}
\end{equation}
Now, choose $\epsilon < \frac{1}{16}$ to conclude
\begin{equation} \label{eq27}
\begin{array} {l}
\frac{d}{dt} \|\sqrt{\rho} \textbf{u}_t (\cdot , t) \|^2 + \widetilde{C }\|
\nabla \textbf{u}_t (\cdot , t) \|^2 \leq
C + C  \| P \Delta \textbf{u} (\cdot , t)\|^2 +
C \|\textbf{u}_t(\cdot , t)\|^2 ,
 \end{array}
\end{equation}
where $\widetilde{C} > 0$ is an absolute constant, and the constant $C$
depends only on $\Omega$, $\| \rho \|_{L^\infty}$, $\| \nabla
\textbf{f} \|$, $\| \textbf{f}_t \|$, $\sup_{t\geq 0} \| \nabla
\textbf{u} \|$. Now, multiplying inequality (\ref{eq27}) by $e^t$
and integrating over $[0 , t ]$, one gets
\begin{eqnarray*}
e^t \|\sqrt{\rho} \textbf{u}_t (\cdot , t)  \|^2 + \widetilde{C } \int_0^t
e^\tau \| \nabla \textbf{u}_t (\cdot , \tau )\|^2 d \tau & \leq &
\|(\sqrt{\rho} \textbf{u}_t) (\cdot , 0)   \|^2 + C  \int_0^t e^\tau \| P \Delta \textbf{u}(\cdot , \tau) \|^2 d \tau
\\
& & \mbox{}  + C  \int_0^t e^\tau d \tau
  + C \int_0^t e^\tau \|\textbf{u}_t (\cdot , \tau) \|^2d \tau  .
\end{eqnarray*}
Hence,
\begin{eqnarray*}
 \|\sqrt{\rho} \textbf{u}_t (\cdot , t) \|^2 +
\widetilde{C } e^{-t}\int_0^t e^\tau \| \nabla \textbf{u}_t (\cdot , \tau) \|^2 d
\tau
& \leq & e^{-t} \beta \| \textbf{u}_t (\cdot , 0)  \|^2
+ C  e^{-t}\int_0^t e^\tau \| P \Delta \textbf{u}(\cdot , \tau) \|^2
d \tau \\
& & \mbox{} + C  e^{-t}\int_0^t e^\tau d \tau
+ C e^{-t} \int_0^t e^\tau \|\textbf{u}_t (\cdot , \tau)\|^2d \tau  .
\end{eqnarray*}
Using inequalities (\ref{eq18}) and (\ref{eq19}), we get the
desired result.
\begin{corollary}\label{cor4} For all $t_0 , t \in \mathbb{R}$, $0 \leq t_0 \leq t$, one has
\begin{eqnarray}
     \int_{t_0}^t \|  \nabla \textbf{u}_t (\cdot , \tau) \|^2 d\tau
  & \leq &  C ( t - t_0 ) + C  . \label{eq30}
\end{eqnarray}
\end{corollary}
\begin{proof} Integrating inequality (\ref{eq27}) from $t_0 $ to $t$, one gets
\begin{eqnarray*}
 \|\sqrt{\rho} \textbf{u}_t (\cdot , t )  \|^2 +
\widetilde{C } \int_{t_0}^t  \| \nabla \textbf{u}_t (\cdot , \tau)\|^2 d \tau & \leq
& \|\sqrt{\rho} \textbf{u}_t  (\cdot , t_0 ) \|^2 + C \int_{t_0}^t d \tau \\
 & & \mbox{} + C  \int_{t_0}^t \| P \Delta \textbf{u}(\cdot , \tau )\|^2 d \tau
 + C \int_{t_0}^t \|\textbf{u}_t (\cdot , \tau) \|^2 d \tau  \\
 & \leq & \|\sqrt{\rho} \textbf{u}_t  (\cdot , t_0 ) \|^2 + C (t - t_0) \\ & &
 \mbox{} +
C  (t - t_0) \left( \sup_{t\geq 0} \| P \Delta \textbf{u}(\cdot , t) \|^2 \right) \\
 & & \mbox{} + C (t - t_0)\left( \sup_{t\geq 0} \|\textbf{u}_t (\cdot , t ) \|^2 \right) . \\
\end{eqnarray*}
Using inequalities (\ref{eq22}) and (\ref{eq23}), one obtains the bound (\ref{eq30}).
\end{proof}
 {\it A priori} estimates for the
solution $\mbox{\boldmath $\xi$}$ of problem (\ref{eq5}), similar to those in Lemma
\ref{lemma4} for $\textbf{u}$, also hold. Indeed, if $\| \nabla
\mbox{\boldmath $\xi$}_0 \| < \delta$, where $\delta$ is the number referred to in
Definition \ref{def1}, then it follows by (\ref{eq7}) that $\|
\nabla \mbox{\boldmath $\xi$} ( \cdot , t )\| \leq \delta M_2$. Therefore, the function
$\widehat{\textbf{u}} = \textbf{u} + \mbox{\boldmath $\xi$}$ is a solution of the
nonhomogeneous Navier-Stokes equations satisfying $\| \nabla
\widehat{\textbf{u}} \| \leq M + \delta M_2$. Moreover, if $\|
P\Delta \mbox{\boldmath $\xi$} ( \cdot , t_0) \|$ is bounded then $\| P\Delta
\widehat{\textbf{u}}(\cdot , t_0) \|$ is also bounded. In this
case, analogously to the proof of Lemma \ref{lemma4}, one can bound
$\| P\Delta \widehat{\textbf{u}} (\cdot , t )\|$, for $t \geq
t_0$. This bound implies that $\| P\Delta \mbox{\boldmath $\xi$} ( \cdot , t) \|$ is
bounded, for $t \geq t_0$. Summarizing,
\begin{lemma}\label{lemma4.5}
   For perturbations $\mbox{\boldmath $\xi$}$ satisfying
   $\| \nabla\mbox{\boldmath $\xi$}_0 \| <\delta$ and
   $\| P\Delta \mbox{\boldmath $\xi$}_0 \| \leq C_0$, $C_0>0$,
   we have $\| P\Delta \mbox{\boldmath $\xi$} (\cdot , t) \| \leq C$, for all $t \geq t_0$. The constant $C$ depends on
   $\| P\Delta \mbox{\boldmath $\xi$}_0\|$, $C_0$, $\Omega$, the initial data of problem (\ref{eq1})
   and on the norms and constants appearing in (\ref{eq6}) and (\ref{eq7}).
\end{lemma}
  It also holds
\begin{lemma} \label{lemma5}
  The function $\mbox{\boldmath $\xi$}$ satisfies
  \begin{equation}
     \int_{t_0}^t \|  \nabla \mbox{\boldmath $\xi$}_t (\cdot , \tau) \|^2 d\tau
   \leq   C ( t - t_0 ) + C , \label{eq31}
  \end{equation}
  for all $t_0$, $t$, $0 \leq t_0 \leq t$.
\end{lemma}
\begin{proof} Note that
$$
  \int_{t_0}^t \|  \nabla \mbox{\boldmath $\xi$}_t (\cdot , \tau) \|^2 d\tau \leq
  C \left( \int_{t_0}^t \|  \nabla \widehat{\textbf{u}}_t (\cdot , \tau) \|^2 d\tau
           + \int_{t_0}^t \|  \nabla \textbf{u}_t (\cdot , \tau) \|^2 d\tau \right).
$$
 The second term on the right hand side of this inequality is bounded
 (Corollary \ref{cor4}). Therefore, it only remains to bound
 $\int_{t_0}^t \|  \nabla \widehat{\textbf{u}}_t (\cdot , \tau) \|^2 d\tau$.
 This bound follows analogously to the
 bound for $\textbf{u}$. \end{proof}

Now, let $\textbf{u} = \sum_{k=1}^{\infty} C_k (t) \textbf{w}^k (x)$ be
the expression of $\textbf{u}$, the solution of (\ref{eq1}), in
terms of the eigenfunctions of the Stokes operator. Let
$\textbf{v}^n := \sum_{k=1}^{n} C_k (t) \textbf{w}^k (x)$ be its
n-th partial sum, and define ${\mathbf e}^n := \textbf{u} - \textbf{v}^n$ and
$\mbox{\boldmath $\mbox{\boldmath $\psi$}$}^n := \textbf{u}^n - \textbf{v}^n$.
We begin by bounding
${\mathbf e}^n$.
\begin{lemma}\label{lemma3} The bounds
\begin{eqnarray}
     \ds \|  \nabla {\mathbf e}^n (\cdot , t) \|^2
  & \leq  & \frac{C}{\lambda_{n+1}}   \label{eq28}
\end{eqnarray} and \begin{eqnarray}
     \ds \|  {\mathbf e}^n (\cdot , t ) \|^2
  & \leq & \frac{C}{\lambda_{n+1}^2} \label{eq29}
\end{eqnarray}
hold for all $t \geq 0$.
\end{lemma}
\begin{proof}
We bound
\begin{eqnarray}
 \| \nabla {\mathbf e}^n (\cdot , t) \|^2 & = & \| \nabla \sum_{k=n+1}^{\infty} C_k (t) \textbf{w}^k (\cdot) \|^2
  \leq  \frac{1}{\lambda_{n+1}}\| P\Delta  \sum_{k=n+1}^{\infty} C_k (t) \textbf{w}^k (\cdot) \|^2 \nonumber\\
 & \leq & \frac{C}{\lambda_{n+1}} \| P\Delta \textbf{u} (\cdot , t )\|^2 \leq \frac{C}{\lambda_{n+1}}, \label{eq21.1}
\end{eqnarray}
as desired.
Moreover,
\begin{eqnarray}
 \|  {\mathbf e}^n (\cdot , t) \|^2 & = & \| \sum_{k=n+1}^{\infty} C_k (t) \textbf{w}^k (\cdot) \|^2
  \leq  \frac{1}{\lambda_{n+1}}\| \nabla \sum_{k=n+1}^{\infty} C_k (t) \textbf{w}^k (\cdot) \|^2 \nonumber\\
 & \leq & \frac{1}{\lambda_{n+1}} \| \nabla {\mathbf e}^n (\cdot , t)\|^2 \leq \frac{C}{\lambda_{n+1}^2}.
 \label{eq21.2}
\end{eqnarray}\end{proof}
Now, we prove that a suitable bound for $\|\nabla \pn(\cdot , t) \|$
implies in a bound for $\| P\Delta \pn (\cdot , t)\|$.
\begin{lemma}\label{lemma7}
   If for some constant $K>0$ the inequality
   $\| \nabla \mbox{\boldmath $\psi$}^n (\cdot ,t)\|^2 \leq \frac{K}{\lambda_{n+1}}$
   holds on an interval
   $[ 0 , t^* ]$, then there exists $C>0$, independent of $n$, such that
  \begin{equation}
    \| P \Delta \pn (\cdot , t )\|^2 \leq C,
  \end{equation}
  for all $ t \in [ 0 , t^* ] $.
\end{lemma}
 {\bf Proof:} Since $\pn = \un - \textbf{v}^n $ and $\sup_{t\geq 0}\| P\Delta \textbf{v}^n (\cdot , t) \| \leq
 \sup_{t\geq 0} \| P\Delta \textbf{u} (\cdot , t )\| < \infty$, one only needs to bound
$ \| P\Delta \un \|$. To this end, note that
$$
 \| \nabla \pn \|^2 =  \| \nabla \un - \nabla \textbf{v}^n \|^2 \leq
 \frac{K}{\lambda_{n+1}}.
$$
Therefore,
$$
 \| \nabla \un \| \leq \left( \frac{K}{\lambda_{n+1}}\right)^\frac{1}{2} +
 \| \nabla \textbf{v}^n \| \leq \left( \frac{K}{\lambda_{1}}\right)^\frac{1}{2} + M .
$$
 We also have, for $\un$, the bounds
 \begin{eqnarray*}
   & & \displaystyle \int_{t_0}^t \| P \Delta \un (\cdot,\tau)\|^2 d \tau \leq C + C(t - t_0) , \\
    & & \displaystyle \int_{t_0}^t \| \un_t (\cdot, \tau) \|^2 d \tau \leq C + C(t - t_0) ,
  \end{eqnarray*}
 which are analogous to the bounds (\ref{eq14}) and (\ref{eq15}) for $\textbf{u}$, and
 can be proved by analogous arguments. Therefore, by Lemma \ref{lemma11},
 \begin{eqnarray*}
   & & \displaystyle \sup_{t \geq 0} e^{-t} \int_{0}^{t} e^{\tau} \| \un_t (\cdot , \tau) \|^2
   d \tau \leq C , \\
   & & \ds \sup_{t \geq 0} e^{-t} \int_{0}^{t} e^{\tau} \| P \Delta \un (\cdot , \tau)\|^2 d \tau
   \leq C  .
 \end{eqnarray*}
 Using these inequalities, one can show that
\begin{equation}\label{eq37}
\begin{array} {l}
\frac{d}{dt} \|\sqrt{\rho} \un_t  \|^2 + \widetilde{C }\| \nabla \un_t \|^2 \leq
C + C  \| P \Delta \un\|^2 + C \|\un_t\|^2,
 \end{array}
\end{equation}
for all $t \in [0 , t^* ]$. At this point, we need to restrict the time interval, since
the constant $C$ depends also on $ \sup_{t\geq 0} \| \nabla \un \|$, and
we can assure this term to be bounded,
uniformly with respect to $n$, only in the interval $[0 , t^*]$.\\
Using inequality (\ref{eq37}), it follows that
\begin{equation}
 \|  \un_t (\cdot , t) \|
    \leq C .  \label{eq38}
\end{equation}
Finally, inequality (\ref{eq38}) allows one to prove
\begin{equation}
       \|  P\Delta \un  (\cdot , t) \|\leq C  ,    \label{eq39}
  \end{equation}
for all $t \in [0 , t^* ]$. We do not give the details of the
proof, since it is completely analogous to the proof of Lemma
\ref{lemma4}. \end{proof}

Now we bound, for later use, the term $\nabla P_n ( \pn - \mbox{\boldmath $\xi$} ) = \nabla
\pn - \nabla P_n \mbox{\boldmath $\xi$}$. First, note that $\textbf{v}^n$ satisfies
\begin{equation} \label{eq40}
  ( \rho \textbf{u}_t , \mbox{\boldmath $\phi$}^n ) + ( \nabla \textbf{v}^n , \nabla \mbox{\boldmath $\phi$}^n )
  + ( \rho \textbf{u} \cdot \nabla \textbf{u} , \mbox{\boldmath $\phi$}^n ) = ( \rho \textbf{f} , \mbox{\boldmath $\phi$}^n ) ,
\end{equation}
for all $\mbox{\boldmath $\phi$}^n$ of the form $\mbox{\boldmath $\phi$}^n ( x) = \sum_{k=1}^{n}
d_k \textbf{w}^k (x)$. Subtracting equation (\ref{eq8}) from
equation (\ref{eq40}), we get
\begin{eqnarray}
  ( \rho^n \pn_t , \mbox{\boldmath $\phi$}^n ) + ( \nabla \pn , \nabla \mbox{\boldmath $\phi$}^n ) & = &
   ( \rho^n {\mathbf e}^n_t , \mbox{\boldmath $\phi$}^n ) + ( (\rho - \rn ) ( \textbf{u}_t + \textbf{u} \cdot \nabla \textbf{u} - \textbf{f}) , \mbox{\boldmath $\phi$}^n)
   \nonumber\\
 & & \mbox{}  + ( \rho^n \textbf{u} \cdot \nabla \textbf{u} , \mbox{\boldmath $\phi$}^n )
 - ( \rho^n \textbf{u}^n \cdot \nabla \textbf{u}^n , \mbox{\boldmath $\phi$}^n ). \label{eq41}
\end{eqnarray}
On the other hand, taking the inner product of the fist equation
in (\ref{eq5}) with $\mbox{\boldmath $\phi$}^n$ and integrating by parts, one gets
\begin{eqnarray}
   ( \widehat{\rho} \mbox{\boldmath $\xi$}_t , \mbox{\boldmath $\phi$}^n ) +
  ( \widehat{\rho} \textbf{u} \cdot \nabla \mbox{\boldmath $\xi$} , \mbox{\boldmath $\phi$}^n ) +
  ( \widehat{\rho} \mbox{\boldmath $\xi$} \cdot \nabla \textbf{u} , \mbox{\boldmath $\phi$}^n ) +
  ( \widehat{\rho} \mbox{\boldmath $\xi$} \cdot \nabla \mbox{\boldmath $\xi$} , \mbox{\boldmath $\phi$}^n ) +
  ( \nabla \mbox{\boldmath $\xi$} , \nabla \mbox{\boldmath $\phi$}^n ) = & & \nonumber \\
   ( (\rho - \rn ) ( \textbf{u}_t + \textbf{u} \cdot \nabla \textbf{u} - \textbf{f}) , \mbox{\boldmath $\phi$}^n) .
   \hspace{.5cm} & & \label{eq42}
\end{eqnarray}
Subtract equation (\ref{eq42}) from equation (\ref{eq41}) to get
\begin{eqnarray}
  ( \rho^n \mbox{\boldmath $\theta$}_t , \mbox{\boldmath $\phi$}^n ) + ( \nabla \mbox{\boldmath $\theta$} , \nabla \mbox{\boldmath $\phi$}^n )  & = &
( (\widehat{\rho} - \rn ) ( \textbf{u}_t + \textbf{u} \cdot \nabla
\textbf{u} - \textbf{f}) , \mbox{\boldmath $\phi$}^n)
+ ((\widehat{\rho} - \rn ) \mbox{\boldmath $\xi$}_t , \mbox{\boldmath $\phi$}^n ) \nonumber \\
 & & \hspace{-1cm} \mbox{} + ( \rn {\mathbf e}^n_t , \mbox{\boldmath $\phi$}^n  ) + ( \rho^n \textbf{u} \cdot \nabla \textbf{u} , \mbox{\boldmath $\phi$}^n )
 - ( \rho^n \textbf{u}^n \cdot \nabla \textbf{u}^n , \mbox{\boldmath $\phi$}^n ) \nonumber\\
& & \hspace{-1cm} \mbox{}
+ ( \widehat{\rho} \textbf{u} \cdot \nabla \mbox{\boldmath $\xi$} , \mbox{\boldmath $\phi$}^n )
 +  ( \widehat{\rho} \mbox{\boldmath $\xi$} \cdot \nabla \textbf{u} , \mbox{\boldmath $\phi$}^n ) +
  ( \widehat{\rho} \mbox{\boldmath $\xi$} \cdot \nabla \mbox{\boldmath $\xi$} , \mbox{\boldmath $\phi$}^n ) , \label{eq43}
\end{eqnarray}
where $\mbox{\boldmath $\theta$} := \mbox{\boldmath $\psi$}^n - \mbox{\boldmath $\xi$}$. Now, since
$ P_n \mbox{\boldmath $\theta$}  = P_n ( \pn - \mbox{\boldmath $\xi$} ) = \pn - P_n \mbox{\boldmath $\xi$}$, one has
\begin{eqnarray*}
 ( \rho^n \textbf{u} \cdot \nabla \pn , \mbox{\boldmath $\phi$}^n ) & = & (\rho^n  \textbf{u} \cdot \nabla P_n \mbox{\boldmath $\theta$} , \mbox{\boldmath $\phi$}^n )
 +  (\rho^n  \textbf{u} \cdot \nabla P_n \mbox{\boldmath $\xi$} , \mbox{\boldmath $\phi$}^n ), \\
 ( \rho^n \pn \cdot \nabla \textbf{u} , \mbox{\boldmath $\phi$}^n ) & = &
 (\rho^n  P_n \mbox{\boldmath $\theta$} \cdot \nabla \textbf{u}  , \mbox{\boldmath $\phi$}^n )
 + (\rho^n  P_n \mbox{\boldmath $\xi$}  \cdot \nabla \textbf{u} , \mbox{\boldmath $\phi$}^n ), \\
( \rho^n \pn \cdot \nabla \pn , \mbox{\boldmath $\phi$}^n ) & = &
 (\rho^n  \pn  \cdot \nabla P_n \mbox{\boldmath $\theta$} , \mbox{\boldmath $\phi$}^n )
 + (\rho^n  \pn   \cdot \nabla P_n \mbox{\boldmath $\xi$} , \mbox{\boldmath $\phi$}^n ).
\end{eqnarray*}
Therefore,
\begin{eqnarray*}
  & & ( \rho^n \textbf{u} \cdot \nabla \textbf{u} , \mbox{\boldmath $\phi$}^n )
 - ( \rho^n \textbf{u}^n \cdot \nabla \textbf{u}^n , \mbox{\boldmath $\phi$}^n ) =
 ( \rn \pn \cdot \nabla {\mathbf e}^n , \mbox{\boldmath $\phi$}^n ) +
 ( \rn {\mathbf e}^n \cdot \nabla \pn , \mbox{\boldmath $\phi$}^n ) \\
 & & \mbox{} + ( \rn \textbf{u} \cdot \nabla {\mathbf e}^n , \mbox{\boldmath $\phi$}^n )
  + ( \rn {\mathbf e}^n \cdot \nabla \textbf{v}^n , \mbox{\boldmath $\phi$}^n ) - ( \rn \textbf{u} \cdot \nabla P_n \mbox{\boldmath $\theta$} , \mbox{\boldmath $\phi$}^n )
 - ( \rn \textbf{u} \cdot \nabla P_n \mbox{\boldmath $\xi$} , \mbox{\boldmath $\phi$}^n ) \\
 & & \mbox{} - ( \rn P_n \mbox{\boldmath $\theta$} \cdot \nabla \textbf{u} , \mbox{\boldmath $\phi$}^n )
- ( \rn P_n \mbox{\boldmath $\xi$} \cdot \nabla \textbf{u} , \mbox{\boldmath $\phi$}^n ) -
( \rn \pn \cdot \nabla P_n \mbox{\boldmath $\theta$} , \mbox{\boldmath $\phi$}^n )
    - ( \rn \pn \cdot \nabla P_n \mbox{\boldmath $\xi$} , \mbox{\boldmath $\phi$}^n ) .
\end{eqnarray*}
Moreover, since $ \mbox{\boldmath $\xi$} = P_n \mbox{\boldmath $\xi$} + Q_n \mbox{\boldmath $\xi$}$,
one can show, after some computations, that
\begin{eqnarray*}
   & & ( \rho^n \textbf{u} \cdot \nabla \textbf{u} , \mbox{\boldmath $\phi$}^n ) - ( \rho^n \textbf{u}^n \cdot \nabla \textbf{u}^n , \mbox{\boldmath $\phi$}^n )
   + ( \widehat{\rho} \textbf{u} \cdot \nabla \mbox{\boldmath $\xi$} , \mbox{\boldmath $\phi$}^n )
   +  ( \widehat{\rho} \mbox{\boldmath $\xi$} \cdot \nabla \textbf{u} , \mbox{\boldmath $\phi$}^n )
    + ( \widehat{\rho} \mbox{\boldmath $\xi$} \cdot \nabla \mbox{\boldmath $\xi$} , \mbox{\boldmath $\phi$}^n )  =  \\
   & & ( \rn \pn \cdot \nabla {\mathbf e}^n , \mbox{\boldmath $\phi$}^n ) +
       ( \rn {\mathbf e}^n \cdot \nabla \pn , \mbox{\boldmath $\phi$}^n ) + ( \rn \textbf{u} \cdot \nabla {\mathbf e}^n , \mbox{\boldmath $\phi$}^n )
       +  ( \rn {\mathbf e}^n \cdot \nabla \textbf{v}^n , \mbox{\boldmath $\phi$}^n ) \\
   & & \mbox{} - ( \rn \textbf{u} \cdot \nabla P_n \mbox{\boldmath $\theta$} , \mbox{\boldmath $\phi$}^n ) - ( \rn P_n \mbox{\boldmath $\theta$} \cdot \nabla \textbf{u} , \mbox{\boldmath $\phi$}^n )
   + ( (\widehat{\rho} - \rn ) \textbf{u} \cdot \nabla P_n \mbox{\boldmath $\xi$} , \mbox{\boldmath $\phi$}^n) \\
   & & \mbox{} + ( (\widehat{\rho} - \rn ) P_n \mbox{\boldmath $\xi$} \cdot \nabla \textbf{u} , \mbox{\boldmath $\phi$}^n)
   + ( \widehat{\rho} \textbf{u} \cdot \nabla Q_n \mbox{\boldmath $\xi$} , \mbox{\boldmath $\phi$}^n) + ( \widehat{\rho} Q_n \mbox{\boldmath $\xi$} \cdot \nabla \textbf{u} , \mbox{\boldmath $\phi$}^n)
   - ( \rn \pn \cdot \nabla P_n \mbox{\boldmath $\theta$} , \mbox{\boldmath $\phi$}^n ) \\
   & & \mbox{}
   - ( \rn P_n \mbox{\boldmath $\theta$} \cdot \nabla P_n \mbox{\boldmath $\xi$} , \mbox{\boldmath $\phi$}^n )
   + ( \rn P_n \mbox{\boldmath $\xi$} \cdot \nabla Q_n \mbox{\boldmath $\xi$} , \mbox{\boldmath $\phi$}^n )
   + ( \rn Q_n \mbox{\boldmath $\xi$} \cdot \nabla P_n \mbox{\boldmath $\xi$} , \mbox{\boldmath $\phi$}^n ) \\
   & & \mbox{}
   +  ( \rn Q_n \mbox{\boldmath $\xi$} \cdot \nabla Q_n \mbox{\boldmath $\xi$} , \mbox{\boldmath $\phi$}^n )
   + ( (\widehat{\rho} - \rn ) \mbox{\boldmath $\xi$} \cdot \nabla \mbox{\boldmath $\xi$} , \mbox{\boldmath $\phi$}^n ).
\end{eqnarray*}
Applying this identity to (\ref{eq43}), and taking
$\mbox{\boldmath $\phi$}^n = P_n \mbox{\boldmath $\theta$}_t $,
one obtains
\begin{eqnarray}
 \| \sqrt{\rn} \mbox{\boldmath $\theta$}_t \|^2 +
  \frac{1}{2} \frac{d}{dt} \| \nabla P_{n}\mbox{\boldmath $\theta$} \|^2  & = &
 ( (\widehat{\rho} - \rn ) ( \textbf{u}_t + \textbf{u} \cdot \nabla \textbf{u}
  - \textbf{f} + \mbox{\boldmath $\xi$}_t
    + \textbf{u} \cdot \nabla P_n \mbox{\boldmath $\xi$} ) , P_n \mbox{\boldmath $\theta$}_t )\nonumber  \label{equacaotheta}
  \\
  & &\mbox{}  + ( (\widehat{\rho} - \rn ) (  P_n \mbox{\boldmath $\xi$} \cdot \nabla \textbf{u} + \mbox{\boldmath $\xi$} \cdot \nabla \mbox{\boldmath $\xi$}) ,
       P_n \mbox{\boldmath $\theta$}_t ) - ( \rn \mbox{\boldmath $\theta$}_t , Q_n \mbox{\boldmath $\xi$}_t  )
        \\
  & &  \mbox{}   + ( \rn {\mathbf e}^n_t , P_n \mbox{\boldmath $\theta$}_t  ) +
  ( \rn \pn \cdot \nabla {\mathbf e}^n , P_n \mbox{\boldmath $\theta$}_t ) +
       ( \rn {\mathbf e}^n \cdot \nabla \pn , P_n \mbox{\boldmath $\theta$}_t ) \nonumber\\
        & &  \mbox{}
        + ( \rn \textbf{u} \cdot \nabla {\mathbf e}^n , P_n \mbox{\boldmath $\theta$}_t ) +
        ( \rn {\mathbf e}^n \cdot \nabla \textbf{v}^n , P_n \mbox{\boldmath $\theta$}_t )
   - ( \rn \textbf{u} \cdot \nabla P_n \mbox{\boldmath $\theta$} , P_n \mbox{\boldmath $\theta$}_t )\nonumber\\
   & & \mbox{}
   - ( \rn P_n \mbox{\boldmath $\theta$} \cdot \nabla \textbf{u} , P_n \mbox{\boldmath $\theta$}_t )
   + ( \widehat{\rho} \textbf{u} \cdot \nabla Q_n \mbox{\boldmath $\xi$} , P_n \mbox{\boldmath $\theta$}_t)
   + ( \widehat{\rho} Q_n \mbox{\boldmath $\xi$} \cdot \nabla \textbf{u} , P_n \mbox{\boldmath $\theta$}_t) \nonumber \\
  & & \mbox{}  - ( \rn \pn \cdot \nabla P_n \mbox{\boldmath $\theta$} , P_n \mbox{\boldmath $\theta$}_t )
   - ( \rn P_n \mbox{\boldmath $\theta$} \cdot \nabla P_n \mbox{\boldmath $\xi$} , P_n \mbox{\boldmath $\theta$}_t ) \nonumber\\
   & & \mbox{} + ( \rn P_n \mbox{\boldmath $\xi$} \cdot \nabla Q_n \mbox{\boldmath $\xi$} , P_n \mbox{\boldmath $\theta$}_t ) +
   ( \rn Q_n \mbox{\boldmath $\xi$} \cdot \nabla P_n \mbox{\boldmath $\xi$} , P_n \mbox{\boldmath $\theta$}_t ) \nonumber \\
   & & \mbox{} +  ( \rn Q_n \mbox{\boldmath $\xi$} \cdot \nabla Q_n \mbox{\boldmath $\xi$} , P_n \mbox{\boldmath $\theta$}_t ).
   \nonumber
\end{eqnarray}
Now, we bound each term on the right hand side of identity (\ref{equacaotheta}). Given $\epsilon >0$, we bound
\begin{eqnarray*}
 | ( \rn \mbox{\boldmath $\theta$}_t , Q_n \mbox{\boldmath $\xi$}_t  ) |
 & \leq &
    \epsilon \| \mbox{\boldmath $\theta$}_t \|^2 + \frac{C(\epsilon)}{\lambda_{n+1}} \| \nabla \mbox{\boldmath $\xi$}_t \|^2, \\
 | ( \rn {\mathbf e}^n_t , P_n \mbox{\boldmath $\theta$}_t  ) |
 & \leq &
    \epsilon \| \mbox{\boldmath $\theta$}_t \|^2 + \frac{C(\epsilon)}{\lambda_{n+1}} \| \nabla \textbf{u}_t \|^2, \\
 | ( \rn \pn \cdot \nabla {\mathbf e}^n , P_n \mbox{\boldmath $\theta$}_t ) |
 & \leq &
    \epsilon \| \mbox{\boldmath $\theta$}_t \|^2 + \frac{C(\epsilon)}{\lambda_{n+1}} \| P\Delta\pn \|^2,\\
 | ( \rn {\mathbf e}^n \cdot \nabla \pn , P_n \mbox{\boldmath $\theta$}_t ) |
 & \leq &
   \epsilon \| \mbox{\boldmath $\theta$}_t \|^2 + \frac{C(\epsilon)}{\lambda_{n+1}} \| P\Delta\pn \|^2,  \\
 | ( \rn \textbf{u} \cdot \nabla {\mathbf e}^n , P_n \mbox{\boldmath $\theta$}_t ) |
 & \leq &
     \epsilon \| \mbox{\boldmath $\theta$}_t \|^2 + \frac{C(\epsilon)}{\lambda_{n+1}},\\
 | ( \rn {\mathbf e}^n \cdot \nabla \textbf{v}^n , P_n \mbox{\boldmath $\theta$}_t ) |
 & \leq &
     \epsilon \| \mbox{\boldmath $\theta$}_t \|^2 + \frac{C(\epsilon)}{\lambda_{n+1}}, \\
 | ( \rn \textbf{u} \cdot \nabla P_n \mbox{\boldmath $\theta$} , P_n \mbox{\boldmath $\theta$}_t ) |
 & \leq &
      \epsilon \| \mbox{\boldmath $\theta$}_t \|^2 + C(\epsilon) \| \nabla P_n\mbox{\boldmath $\theta$} \|^2, \\
 | ( \rn P_n \mbox{\boldmath $\theta$} \cdot \nabla \textbf{u} , P_n \mbox{\boldmath $\theta$}_t ) |
 & \leq &
    \epsilon \| \mbox{\boldmath $\theta$}_t \|^2 + C(\epsilon) \| \nabla P_n\mbox{\boldmath $\theta$} \|^2, \\
 | ( \widehat{\rho} \textbf{u} \cdot \nabla Q_n \mbox{\boldmath $\xi$} , P_n \mbox{\boldmath $\theta$}_t) |
 & \leq &
    \epsilon \| \mbox{\boldmath $\theta$}_t \|^2 + \frac{C(\epsilon)}{\lambda_{n+1}} \| P\Delta \mbox{\boldmath $\xi$} \|^2, \\
 | ( \widehat{\rho} Q_n \mbox{\boldmath $\xi$} \cdot \nabla \textbf{u} , P_n \mbox{\boldmath $\theta$}_t) |
 & \leq &
    \epsilon \| \mbox{\boldmath $\theta$}_t \|^2 + \frac{C(\epsilon)}{\lambda_{n+1}} \| P\Delta \mbox{\boldmath $\xi$} \|^2, \\
 | ( \rn \pn \cdot \nabla P_n \mbox{\boldmath $\theta$} , P_n \mbox{\boldmath $\theta$}_t ) |
 & \leq &
   \epsilon \| \mbox{\boldmath $\theta$}_t \|^2 + C(\epsilon) \| P\Delta \pn  \|^2 \| \nabla P_n \mbox{\boldmath $\theta$} \|^2,\\
 | ( \rn P_n \mbox{\boldmath $\theta$} \cdot \nabla P_n \mbox{\boldmath $\xi$} , P_n \mbox{\boldmath $\theta$}_t ) |
 & \leq &
   \epsilon \| \mbox{\boldmath $\theta$}_t \|^2 + C(\epsilon) \| P\Delta \mbox{\boldmath $\xi$}  \|^2 \| \nabla P_n \mbox{\boldmath $\theta$}  \|^2, \\
 |( \rn P_n \mbox{\boldmath $\xi$} \cdot \nabla Q_n \mbox{\boldmath $\xi$} , P_n \mbox{\boldmath $\theta$}_t ) |
 & \leq &
   \epsilon \| \mbox{\boldmath $\theta$}_t \|^2 + \frac{C(\epsilon)}{\lambda_{n+1}} \| P\Delta \mbox{\boldmath $\xi$} \|^4, \\
 |( \rn Q_n \mbox{\boldmath $\xi$} \cdot \nabla P_n \mbox{\boldmath $\xi$} , P_n \mbox{\boldmath $\theta$}_t ) |
 & \leq &
   \epsilon \| \mbox{\boldmath $\theta$}_t \|^2 + \frac{C(\epsilon)}{\lambda_{n+1}} \| P\Delta \mbox{\boldmath $\xi$} \|^4, \\
 | ( \rn Q_n \mbox{\boldmath $\xi$} \cdot \nabla Q_n \mbox{\boldmath $\xi$} , P_n \mbox{\boldmath $\theta$}_t ) |
 & \leq &
   \epsilon \| \mbox{\boldmath $\theta$}_t \|^2 + \frac{C(\epsilon)}{\lambda_{n+1}} \| P\Delta \mbox{\boldmath $\xi$} \|^4 .
\end{eqnarray*}
 It remains to bound $| ( \pi g^n , P_n \mbox{\boldmath $\theta$}_t ) |$, where
 $ \pi := \widehat{\rho} - \rn$ and \\ ${\mathbf g}^n := \textbf{u}_t + \textbf{u} \cdot \nabla \textbf{u} - \textbf{f} + \mbox{\boldmath $\xi$}_t
    + \textbf{u} \cdot \nabla P_n \mbox{\boldmath $\xi$} + P_n \mbox{\boldmath $\xi$} \cdot \nabla \textbf{u} + \mbox{\boldmath $\xi$} \cdot \nabla \mbox{\boldmath $\xi$}$. We begin by
    bounding ${\mathbf g}^n$.
\begin{lemma}\label{lemma7.5}
  For all $p$, $2\leq p \leq 6$, the bound
  \begin{eqnarray}
   \| {\mathbf g}^n (\cdot, t) \|_{L^p}^2 & \leq & C + C \| P\Delta \mbox{\boldmath $\xi$}(\cdot , t) \|^2
   + C \| P\Delta \mbox{\boldmath $\xi$} (\cdot , t) \|^4 \nonumber \\
   & & \mbox{}+ C \| \nabla \textbf{u}_t (\cdot , t ) \|^{2} + C \| \nabla \mbox{\boldmath $\xi$}_t (\cdot , t ) \|^{2}\label{eq44}
   \end{eqnarray}
  holds for all $t\geq 0$.
\end{lemma}
\begin{proof}Since $2\leq p\leq 6$, we have
  \begin{eqnarray*}
  \| {\mathbf g}^n (\cdot, t) \|_{L^p}^2
  & \leq &
  C \left\{ \| \textbf{u}_t (\cdot, t) \|_{L^p}^2 + \| \textbf{u} \cdot \nabla \textbf{u}(\cdot, t)\|_{L^p}^2
 +\| \textbf{f}(\cdot, t)\|_{L^p}^2 + \| \mbox{\boldmath $\xi$}_t (\cdot, t)\|_{L^p}^2\right. \\
 & & \left. \hspace{.6cm}\mbox{} + \| \textbf{u} \cdot \nabla P_n
 \mbox{\boldmath $\xi$}(\cdot, t)\|_{L^p}^2
 + \| P_n \mbox{\boldmath $\xi$} \cdot \nabla \textbf{u}(\cdot, t)\|_{L^p}^2
 + \| \mbox{\boldmath $\xi$} \cdot \nabla \mbox{\boldmath $\xi$} (\cdot, t)\|_{L^p}^2 \right\} \\
  & \leq  & C\left\{\| \nabla \textbf{u}_t (\cdot, t)\|^2 +\| P \Delta \textbf{u} (\cdot, t)\|^4 +
  \| \nabla \textbf{f}(\cdot, t) \|^2 + \| \nabla \mbox{\boldmath $\xi$}_t (\cdot, t)\|^2 \right.\\
  & & \left. \hspace{.6cm}\mbox{} + \| P \Delta \textbf{u} (\cdot , t )\|^2
  \| P \Delta \mbox{\boldmath $\xi$} (\cdot , t )\|^2
  + \| P \Delta \mbox{\boldmath $\xi$} (\cdot , t ) \|^4 \right\} \\
  & \leq &  C + C \| P\Delta \mbox{\boldmath $\xi$} (\cdot , t) \|^2
   + C \| P\Delta \mbox{\boldmath $\xi$} (\cdot , t) \|^4
   + C \| \nabla \textbf{u}_t (\cdot , t)\|^2 + C\| \nabla \mbox{\boldmath $\xi$}_t (\cdot , t)\|^2 .
\end{eqnarray*}
\end{proof}
\begin{lemma}\label{lemma8}
 If $6 \leq p_0 < \infty$, then the bound
 \begin{eqnarray}
   \| \pi (\cdot , t) \|_{L^r}^2 & \leq & C \|\pi (\cdot , t_0 ) \|_{L^r}^2
   + C (t-t_0 ) \int_{t_0}^{t} \| \nabla P_n \mbox{\boldmath $\theta$} (\cdot , \tau )\|^2 d \tau \nonumber \\
  & & \hspace{-.5cm}\mbox{} + \frac{C}{\lambda_{n+1}} (t-t_0 )^2
  + \frac{C}{\lambda_{n+1}} (t-t_0 ) \int_{t_0}^{t} \| P \Delta \mbox{\boldmath $\xi$}(\cdot , \tau )\|^2 d \tau \label{eq45}
 \end{eqnarray}
 holds for all $t\geq 0$ and all $r$, $2\leq r \leq \frac{6p_0}{6 + p_0}$. If $p_0 = \infty$, then
 the bound (\ref{eq45}) is valid for all $t\geq 0$ and all  $r$, $2\leq r \leq 6$.
 \end{lemma}
\begin{proof}
 First note that
 $$
  \pi_t + \textbf{u}^n \cdot \nabla \pi = ( \textbf{u}^n - \widehat{\textbf{u}}) \cdot \nabla \widehat{\rho}.
 $$
 Since $\widehat{\textbf{u}}  = \textbf{u} + \mbox{\boldmath $\xi$}$, we write the equation above as
 \begin{equation} \label{eq59}
  \pi_t + \textbf{u}^n \cdot \nabla \pi = (P_n \mbox{\boldmath $\theta$} - Q_n \mbox{\boldmath $\xi$} - {\mathbf e}^n) \cdot \nabla \widehat{\rho}.
 \end{equation}
 Let $r$ belonging to the suitable interval depending on the value of $p_0$.
 Multiply equation (\ref{eq59}) by $| \pi |^{r-1}$ and integrate to get
 \begin{eqnarray*}
   \frac{1}{r} \frac{d}{dt} \| \pi (\cdot , t ) \|_{L^r}^r & \leq &
   \int_{\Omega} (P_n \mbox{\boldmath $\theta$} - Q_n \mbox{\boldmath $\xi$} -
   {\mathbf e}^n) \cdot \nabla \widehat{\rho}| \pi |^{r-1} dx \\
   & \leq & \left( \int_{\Omega} | (P_n \mbox{\boldmath $\theta$} - Q_n \mbox{\boldmath $\xi$} -
   {\mathbf e}^n)\cdot \nabla \widehat{\rho}|^r dx \right)^{\frac{1}{r}}
   \left( \int_{\Omega} | \pi |^r dx \right)^{\frac{r-1}{r}} \\
   & \leq & \| (P_n \mbox{\boldmath $\theta$} - Q_n \mbox{\boldmath $\xi$}
   - {\mathbf e}^n) \cdot \nabla \widehat{\rho} \|_{L^r} \| \pi \|_{L^r}^{r-1}.
 \end{eqnarray*}
Therefore,
 \begin{eqnarray*}
   \frac{d}{dt} \| \pi (\cdot , t ) \|_{L^r} & \leq &
   \| (P_n \mbox{\boldmath $\theta$} - Q_n \mbox{\boldmath $\xi$} - {\mathbf e}^n) \cdot \nabla \widehat{\rho} \|_{L^r}
   \leq
   \| P_n \mbox{\boldmath $\theta$} - Q_n \mbox{\boldmath $\xi$} - {\mathbf e}^n\|_{L^p}
   \| \nabla \widehat{\rho} \|_{L^{p_0}},
 \end{eqnarray*}
 where $p$ is chosen as $\frac{1}{p} = \frac{1}{r} - \frac{1}{p_0}$ if $6\leq p_0 <\infty$, and as
 $p = r$ if $p_0 = \infty$. Note that in the case $6\leq p_0 <\infty$,
 this choice of $p$ implies $2< \frac{2p_0}{p_0 -2} \leq p \leq 6$. In the case $p_0 = \infty$,
 we have $2 \leq p \leq 6$. In both cases, $p \in [2 , 6]$ and we bound
 \begin{eqnarray*}
  \frac{d}{dt} \| \pi (\cdot , t ) \|_{L^r} & \leq & C \left( \|\nabla P_n \mbox{\boldmath $\theta$}(\cdot , t ) \|
   + \| \nabla Q_n \mbox{\boldmath $\xi$}(\cdot , t )\| +
   \|\nabla {\mathbf e}^n (\cdot , t )\| \right)  \\
   & \leq &
   C
   \left( \|\nabla P_n \mbox{\boldmath $\theta$} (\cdot , t )\| + \frac{C}{(\lambda_{n+1})^{\frac{1}{2}}} +
   \frac{C}{(\lambda_{n+1})^{\frac{1}{2}}} \| P \Delta \mbox{\boldmath $\xi$}(\cdot , t ) \|  \right),
    \end{eqnarray*}
where, for the last inequality, we used Lemma \ref{lema_rautmann}
and inequality (\ref{eq28}). Integrating this inequality from
$t_0$ to $t$,
  \begin{eqnarray*}
   \| \pi (\cdot , t ) \|_{L^r}^2 & \leq & C \left\{ \| \pi (\cdot , t_0 ) \|_{L^r}^2
   + \left( \int_{t_0}^{t} \|\nabla P_n \mbox{\boldmath $\theta$} (\cdot , \tau) \| d\tau \right)^2
   + \left( \int_{t_0}^{t} \frac{1}{(\lambda_{n+1})^{\frac{1}{2}}} d\tau \right)^2 \right. \\
   & & \hspace{.6cm}\left. \mbox{} + \left( \int_{t_0}^{t}  \frac{1}{(\lambda_{n+1})^{\frac{1}{2}}}
      \| P \Delta \mbox{\boldmath $\xi$} (\cdot , \tau )\|d\tau  \right)^2 \right\} \\
   & \leq &  C \left\{ \| \pi (\cdot , t_0 ) \|_{L^r}^2
             + (t - t_0 ) \int_{t_0}^{t} \|\nabla P_n \mbox{\boldmath $\theta$} (\cdot , \tau) \|^2 d\tau
             + \frac{1}{\lambda_{n+1}} (t - t_0 )^2 \right.\\
   & & \hspace{.6cm} \left.
 + \frac{1}{\lambda_{n+1}} (t - t_0 ) \int_{t_0}^{t} \|P\Delta\mbox{\boldmath $\xi$} (\cdot , \tau) \|^2 d\tau \right\},
    \end{eqnarray*}
 which is the desired bound. \end{proof}
 Getting back to inequality (\ref{equacaotheta}), we have
\begin{eqnarray*}
 \alpha \|  \mbox{\boldmath $\theta$}_t \|^2 + \frac{1}{2} \frac{d}{dt} \| \nabla P_n \mbox{\boldmath $\theta$} \|^2 & \leq &
 15 \epsilon \| \mbox{\boldmath $\theta$}_t \|^2 + \frac{C(\epsilon)}{\lambda_{n+1}} \| \nabla \mbox{\boldmath $\xi$}_t \|^2
 + \frac{C(\epsilon)}{\lambda_{n+1}} \| \nabla \textbf{u}_t \|^2
 + \frac{C(\epsilon)}{\lambda_{n+1}} \| P\Delta\pn \|^2\\
 & &
     \mbox{} + \frac{C(\epsilon)}{\lambda_{n+1}}
     + \frac{C(\epsilon)}{\lambda_{n+1}} \| P\Delta \mbox{\boldmath $\xi$} \|^2
     + \frac{C(\epsilon)}{\lambda_{n+1}} \| P\Delta \mbox{\boldmath $\xi$} \|^4
     + C(\epsilon) \| \pi \|_{L^r}^2 \| {\mathbf g}^n \|_{L^p}^2 \\
 & &
 \mbox{} + C(\epsilon) \| \nabla P_n \mbox{\boldmath $\theta$}  \|^2
   \left\{ 1 +\| P\Delta \pn  \|^2 + \| P\Delta \mbox{\boldmath $\xi$}  \|^2 \right\},
\end{eqnarray*}
where, in the case $6 \leq p_0 <\infty$,
the inequality above holds for each $r \in \left[ 3 , \frac{6p_0}{6+p_0}\right]$, with
$p \in \left[\frac{3p_0}{p_0 - 3} , 6\right]$ chosen
such that $\frac{1}{r} + \frac{1}{p} = \frac{1}{2}$. In the case $p_0 = \infty$,
it holds for all $r \in [3 , 6]$, with $p \in [3 , 6]$ chosen such that
$\frac{1}{r} + \frac{1}{p} = \frac{1}{2}$.
 Now fix $\epsilon = \frac{1}{15} \left( \alpha - \frac{1}{2} \right)$. Integrating the
 inequality from $t_0$ to $t$, and using Lemma \ref{lemma4.5},
 Lemma \ref{lemma7} and inequality (\ref{eq31}), we get
\begin{eqnarray}\label{eq46}
 & & \| \nabla P_n \mbox{\boldmath $\theta$} (\cdot , t )\|^2  +
 \int_{t_0}^t \|  \mbox{\boldmath $\theta$}_t ( \cdot , \tau ) \|^2 d \tau \leq
 \| \nabla P_n \mbox{\boldmath $\theta$} (\cdot , t_0 )\|^2 + \frac{C}{\lambda_{n+1}} ( t-t_0 )
 + \frac{C}{\lambda_{n+1}} \int_{t_0}^t \| \nabla \mbox{\boldmath $\xi$}_t (\cdot , \tau ) \|^2 d\tau \nonumber \\
 & & \mbox{} + \frac{C}{\lambda_{n+1}} \int_{t_0}^t \| \nabla \textbf{u}_t (\cdot , \tau ) d \tau\|^2
     + C \int_{t_0}^t  \| \pi (\cdot , \tau ) \|_{L^r}^2 \| {\mathbf g}^n (\cdot , \tau )\|_{L^p}^2 d\tau
     + C \int_{t_0}^t \| \nabla P_n \mbox{\boldmath $\theta$} (\cdot , \tau ) \|^2 d\tau .
\end{eqnarray}
 Adding inequalities (\ref{eq45}) and (\ref{eq46}), one obtains
 \begin{eqnarray}
 & &
     \| \nabla P_n \mbox{\boldmath $\theta$} (\cdot , t )\|^2 + \| \pi (\cdot , t) \|_{L^r}^2
     + \int_{t_0}^t \|  \mbox{\boldmath $\theta$}_t ( \cdot , \tau ) \|^2 d \tau\leq
     \| \nabla P_n \mbox{\boldmath $\theta$} (\cdot , t_0 )\|^2 + C \|\pi (\cdot , t_0 ) \|_{L^r}^2  \label{eq47} \\
 & &
    \mbox{}  + \frac{C}{\lambda_{n+1}} ( t-t_0 ) + \frac{C}{\lambda_{n+1}} (t-t_0 )^2
     + \frac{C}{\lambda_{n+1}} \int_{t_0}^t \| \nabla \mbox{\boldmath $\xi$}_t (\cdot , \tau ) \|^2 d\tau
     + \frac{C}{\lambda_{n+1}} \int_{t_0}^t \| \nabla \textbf{u}_t (\cdot , \tau ) d \tau\|^2 \nonumber \\
 & &
    \mbox{} + C \int_{t_0}^t  \| \pi (\cdot , \tau ) \|_{L^r}^2 \| {\mathbf g}^n (\cdot , \tau )\|_{L^p}^2 d\tau
     + C \int_{t_0}^t \| \nabla P_n \mbox{\boldmath $\theta$} (\cdot , \tau ) \|^2 d\tau
    + C (t-t_0 ) \int_{t_0}^{t} \| \nabla P_n \mbox{\boldmath $\theta$} (\cdot , \tau )\|^2 d \tau . \nonumber
 \end{eqnarray}
Fixing $\bar{t} > t_0$, and using Corollary \ref{cor4} and Lemma \ref{lemma5}, one concludes that
 \begin{eqnarray}
 & &
      \| \nabla P_n \mbox{\boldmath $\theta$} (\cdot , t )\|^2 + \| \pi (\cdot , t) \|_{L^r}^2
     + \int_{t_0}^t \|  \mbox{\boldmath $\theta$}_t ( \cdot , \tau ) \|^2 d \tau
     \leq \| \nabla P_n \mbox{\boldmath $\theta$} (\cdot , t_0 )\|^2 \nonumber\\
 & &
    \mbox{} + C \|\pi (\cdot , t_0 ) \|_{L^r}^2
    + \frac{C}{\lambda_{n+1}} + \frac{C}{\lambda_{n+1}} ( t-t_0 )
    + \frac{C}{\lambda_{n+1}} (t-t_0 )^2 \label{eq48} \\
 & &
    \mbox{} + C \int_{t_0}^t  \| \pi (\cdot , \tau ) \|_{L^r}^2 \| {\mathbf g}^n (\cdot , \tau )\|_{L^p}^2 d\tau
     + C \int_{t_0}^t \{ 1 + (\bar{t} - t_0 ) \} \| \nabla P_n \mbox{\boldmath $\theta$} (\cdot , \tau ) \|^2 d\tau, \nonumber
 \end{eqnarray}
for all $t \in [t_0 , \bar{t} \,]$.
Let $\Lambda ( t) := \| \nabla P_n \mbox{\boldmath $\theta$} (\cdot , t )\|^2 + \| \pi (\cdot , t) \|_{L^r}^2
     +\int_{t_0}^t \|  \mbox{\boldmath $\theta$}_t ( \cdot , \tau ) \|^2 d \tau$.
Thus, inequality (\ref{eq48}) gives
\begin{eqnarray*}
& &
     \Lambda(t)  \leq C \Lambda( t_0 ) +
      \frac{C}{\lambda_{n+1}} + \frac{C}{\lambda_{n+1}} ( t-t_0 )
    + \frac{C}{\lambda_{n+1}} (t-t_0 )^2 \\
& &
    \mbox{} + C \int_{t_0}^t  \left\{ 1 + \bar{t} - t_0 + \| {\mathbf g}^n (\cdot , \tau )\|_{L^p}^2\right\}
 \Lambda (\tau ) d\tau .
\end{eqnarray*}
 Applying a
 corollary of Gronwall's Lemma(see \cite{amann}, page 90, corollary 6.2),
 one gets
 $$
  \Lambda (t) \leq \left( C \Lambda( t_0 ) +
      \frac{C}{\lambda_{n+1}} + \frac{C}{\lambda_{n+1}} ( t-t_0 )
    + \frac{C}{\lambda_{n+1}} (t-t_0 )^2 \right)
    \exp\left\{ C \int_{t_0}^t ( 1 + \bar{t} - t_0 + \| g^n (\cdot , \tau )\|_{L^p}^2) d\tau\right\}.
 $$
We summarize the results in the following lemma.
\begin{lemma}\label{lemma9}
Suppose $\|\nabla \pn (\cdot , t) \|\leq K/ \lambda_{k+1}$ holds
for a constant $K>0$ and
all $t$ in a given interval $0 \leq t_0 \leq t \leq \bar{t}$.
Let $\mbox{\boldmath
$\xi$}$ as in problem (\ref{eq5}), and the functions $\pi$,
$\mbox{\boldmath $\theta$}$, ${\mathbf g}^n$ defined as before. If
$6 \leq p_0 <\infty$ then, for all $t \in [t_0 , \bar{t}\, ]$, one
has
\begin{eqnarray}
& &  \| \nabla P_n \mbox{\boldmath $\theta$} (\cdot , t )\|^2 + \| \pi (\cdot , t) \|_{L^r}^2
     + \int_{t_0}^t \|  \mbox{\boldmath $\theta$}_t ( \cdot , \tau ) \|^2 d \tau \label{eq49}\\
& &
     \leq C \left( \| \nabla P_n \mbox{\boldmath $\theta$} (\cdot , t_0 )\|^2 + \| \pi (\cdot , t_0) \|_{L^r}^2
      + \frac{1}{\lambda_{n+1}} + \frac{1}{\lambda_{n+1}} ( t-t_0 )
    + \frac{1}{\lambda_{n+1}} (t-t_0 )^2 \right) e^{C \int_{t_0}^t a (\tau )d\tau }, \nonumber
\end{eqnarray}
for all $ r \in \left[ 3 , \frac{6p_0}{6+p_0} \right]$, and $p \in
\left[\frac{3p_0}{p_0 - 3} , 6\right]$ chosen such that $\frac{1}{r}
+ \frac{1}{p} = \frac{1}{2}$, where $a(t) := 1 + \bar{t} - t_0 + \|
{\mathbf g}^n (\cdot , t )\|_{L^p}^2$. If $p_0 = \infty$, then the
bound holds for all $r \in [3 , 6]$ and $p\in [3 , 6]$ such that
$\frac{1}{r} + \frac{1}{p} = \frac{1}{2}$.
\end{lemma}
\begin{remark}
From inequalities (\ref{eq44}), (\ref{eq30}), (\ref{eq31}), and Lemma \ref{lemma4.5}, one can bound
$$
 \int_{t_0}^t a(\tau )d\tau \leq \widetilde{C} + \widetilde{C} (t - t_0) + \widetilde{C} (\bar{t} - t_0 ) ( t-t_0 ).
$$
Therefore,
\begin{equation}\label{eq50}
  \exp\left\{ \int_{t_0}^t a(\tau )d\tau \right\}  \leq
     \exp\left\{ \widetilde{C} + \widetilde{C} (t - t_0) + \widetilde{C} (\bar{t} - t_0 ) ( t-t_0 ) \right\}.
\end{equation}
\end{remark}
From now on, we fix the constants $C$ and $\widetilde{C}$
appearing in inequalities (\ref{eq49}) and (\ref{eq50}). We prove now
that the bound for $\| \nabla \pn \|$ required in lemma (\ref{lemma7}) and lemma (\ref{lemma9}) hold for $n$ large enough.
\begin{proposition}\label{claim1}
There exist $K > 0$ and $N \in \mathbb{N}$ such that if $n \geq
N$, then $\|\nabla \pn (\cdot , t )\|^2 < \frac{K}{\lambda_{n+1}}$
for all $t\geq 0$.
\end{proposition}
\begin{proof}
Choose $T$ such that $M_2^2 (F(T))^2 \leq \frac{1}{4}$. Let
$ K := 8 C ( 1 + T + T^2) \exp\{ \widetilde{C} + \widetilde{C} T + \widetilde{C} T^2 \}$ and
let $N$ to be large enough such that $\frac{K}{\lambda_{n+1}} < \delta$ if $n \geq N$. Under these conditions,
we have
\begin{equation}\label{eq51}
 \| \nabla \pn (\cdot , t ) \|^2 < \frac{K}{\lambda_{n+1}},
\end{equation}
for all $t\geq0$. Indeed, suppose that inequality (\ref{eq51}) does not hold. Thus, there exist $n \geq N$ and $t^*>0$ such that
\begin{equation}\label{eq52}
 \| \nabla \mbox{\boldmath $\psi$}^{n} (\cdot , t^*) \|^2 = \frac{K}{\lambda_{n+1}}.
 \end{equation}
 We may have either $t^* \leq T$ or $t^* >T$.
 If $t^* \leq T$, consider $t_0 = 0$, $\mbox{\boldmath $\xi$} = 0$, $\eta = 0$, $\bar{t} = t^*$.
 In this
 case, $\| \nabla P_n \mbox{\boldmath $\theta$} \| = \| \nabla \mbox{\boldmath $\psi$}^{n} \|$.
 Therefore, using Lemma \ref{lemma9}, we have
 \begin{eqnarray*}
  \| \nabla \mbox{\boldmath $\psi$}^{n} (\cdot , t^* )\|^2
+ \| \pi (\cdot , t^*) \|_{L^r}^2
     + \int_{0}^{t^*} \|  \mbox{\boldmath $\psi$}^{n}_t ( \cdot , \tau ) \|^2 d \tau
     & \leq  & \left(  \frac{C}{\lambda_{n+1}} + \frac{C}{\lambda_{n+1}} T
    + \frac{C}{\lambda_{n+1}} T^2 \right) e^{\widetilde{C} + \widetilde{C} T  + \widetilde{C} T^2}, \nonumber \\
    & = & \frac{K}{8 \lambda_{n+1}} < \frac{K}{ \lambda_{n+1}},
\end{eqnarray*}
which contradicts (\ref{eq52}).
If $t^*  > T$, apply Lemma \ref{lemma9}
with $\bar{t} = t^*$, $t_0 = t^* - T$ and $\mbox{\boldmath $\xi$} (x,t)$, $\eta (x,t)$ satisfying
\begin{eqnarray*}
  \mbox{\boldmath $\xi$} (x , t_0 ) & = & \mbox{\boldmath $\psi$}^{n} (x , t_0), \\
  \eta ( x , t_0 ) & = & \rn ( x , t_0 ) - \rho ( x , t_0),
\end{eqnarray*}
to get
\begin{eqnarray*}
& &
\| \nabla \mbox{\boldmath $\psi$}^{n} ( \cdot , t^* ) - \nabla P_n \mbox{\boldmath $\xi$} ( \cdot , t^* ) \|^2
 + \| \rho (\cdot , t^*)- \rho^n (\cdot , t^*)+ \eta (\cdot , t^*) \|_{L^r}^2 \\
& & \mbox{}
     + \int_{t^* - T}^{t^*}
     \| \mbox{\boldmath $\psi$}^{n}_t ( \cdot , \tau ) - \mbox{\boldmath $\xi$}_t ( \cdot , \tau ) \|^2 d \tau
     \leq \left(  \frac{C}{\lambda_{n+1}} + \frac{C}{\lambda_{n+1}} T
    + \frac{C}{\lambda_{n+1}} T^2 \right) e^{\widetilde{C} + \widetilde{C} T
    + \widetilde{C} T^2} = \frac{K}{8 \lambda_{n+1}} .
\end{eqnarray*}
Therefore,
\begin{eqnarray*}
   \| \nabla \pn (\cdot , t^* ) \|^2 & \leq &
   2 \left(\| \nabla \mbox{\boldmath $\psi$}^{n} ( \cdot , t^* ) - \nabla P_n \mbox{\boldmath $\xi$} ( \cdot , t^* ) \|^2
       + \| \nabla P_n \mbox{\boldmath $\xi$} ( \cdot , t^* ) \|^2 \right) \\
   & \leq & 2\left( \frac{K}{8 \lambda_{n+1}} +  M_2^2 \| \nabla \mbox{\boldmath $\xi$} (\cdot , t_0 )\|^2
     F(T)^2 \right) \\
   & \leq & 2\left( \frac{K}{8 \lambda_{n+1}} + \frac{K}{4 \lambda_{n+1}} \right)
    = \frac{3}{4} \frac{K}{ \lambda_{n+1}} < \frac{K}{ \lambda_{n+1}}
\end{eqnarray*}
which again contradicts (\ref{eq52}). This proves the proposition.
\end{proof}
\section{Proof of Theorem \ref{theorem1}} \label{proof}
Using the bounds (\ref{eq28}) and (\ref{eq51}), we directly get
 \begin{equation} \label{eq52.5}
  \| \nabla \textbf{u}(\cdot , t) - \nabla \textbf{u}^n  (\cdot , t) \|^2  \leq
  2 ( \| \nabla \pn (\cdot , t)\|^2 + \| \nabla {\mathbf e}^n (\cdot , t)\|^2  )
   \leq \frac{C}{\lambda_{n+1}},
 \end{equation}
which is the first bound stated in Theorem \ref{theorem1}. In order to prove
the bound (\ref{eq11}) for the density, first note that
\begin{gather}
 \rho_t + \textbf{u} \cdot \nabla \rho  =  0 \label{eq53} \\
 \rho_t^n + \textbf{u}^n \cdot \nabla \rho^n  =  0 .\label{eq54}
\end{gather}
Subtracting equation (\ref{eq54}) from equation (\ref{eq53}), one gets
\begin{equation}\label{eq55}
  (\rho - \rho^n )_t + \textbf{u}^n \cdot \nabla (\rho - \rho^n ) = (\textbf{u}^n - \textbf{u} ) \cdot \nabla \rho.
\end{equation}
Now, if $6\leq p_0 <\infty$, let $r \in \left[2 , \frac{6p_0}{6+p_0} \right]$. Choose
$p \in \left[\frac{2p_0}{p_0 -2} , 6 \right]$ such that $\frac{1}{r} = \frac{1}{p} + \frac{1}{p_0}$.
Multiplying equation (\ref{eq55}) by $ | \rho - \rho^n |^{r-1}$ and integrating
over $\Omega$, we obtain
\begin{eqnarray*}
  \frac{1}{r} \frac{d}{dt} \| \rho - \rho^n \|^r_{L^r} & = &
     \int_{\Omega} | \rho - \rho^n |^{r-1} (\textbf{u}^n - \textbf{u} ) \cdot \nabla \rho dx \\
     & \leq & \| \rho - \rho^n \|_{L^r}^{r-1} \| (\textbf{u}^n - \textbf{u} ) \cdot \nabla \rho\|_{L^r} \\
     & \leq & \| \rho - \rho^n \|^{r-1}_{L^r}\| \nabla \rho \|_{L^{p_0}} \| \textbf{u}^n - \textbf{u} \|_{L^p},
\end{eqnarray*}
If $p_0 = \infty$, the bounds above hold for all $r \in [ 2, 6]$ and $p = r$.
Thus,
\begin{equation}\label{eq56}
 \frac{d}{dt} \| \rho(\cdot , t ) - \rho^n (\cdot , t) \|_{L^r} \leq
 C \| \textbf{u}^n (\cdot ,t) - \textbf{u} (\cdot , t) \|_{L^r} \leq
  C \| \nabla \textbf{u}^n (\cdot , t ) - \nabla \textbf{u} (\cdot , t ) \|.
\end{equation}
Integrating inequality (\ref{eq56}) from $0$ to $t$ and using $(\ref{eq52.5})$, one gets
\begin{equation}\label{eq58}
 \| \rho (\cdot , t )- \rho^n (\cdot , t) \|_{L^r} \leq
  \frac{C}{(\lambda_{n+1})^{\frac{1}{2}}} t
  + \| \rho (\cdot , 0 )- \rho^n (\cdot , 0) \|_{L^r} =
  \frac{C}{(\lambda_{n+1})^{\frac{1}{2}}} t,
\end{equation}
since $\rho^n ( x , 0 ) = \rho_0 (x) $.
This finishes the proof of Theorem \ref{theorem1}.
\appendix

\section{Proof of Lemma \ref{lemma11}} \label{apendice}

Suppose $h(t)$ integrable, nonnegative and satisfying, for all $t$, $t_0$ with $0\leq t_0 \leq t$,
\begin{equation}\label{eq57}
 \displaystyle \int_{t_0}^t h (\tau) d \tau \leq a_1 (t - t_0) + a_2,
\end{equation}
 where $a_1$ and $a_2$ are nonnegative constants. We first consider the case $0\leq t \leq 1$.
 In this case,
 \[
  e^{-t}\int_{t_0}^t e^\tau h (\tau) d \tau \leq \int_{t_0}^t h (\tau) d \tau \leq
  \int_{0}^1 h (\tau) d \tau \leq a_1 + a_2.
 \]
 Now, if $t>1$, let $n \in \mathbb{N}$ and $r \in [0 , 1)$ such that $t = n + r$. Then,
 \begin{eqnarray*}
  \int_0^t e^\tau h(\tau) d\tau & = &  \int_0^{n+r} e^\tau h(\tau) d\tau =
  \sum_{j=1}^n \int_{j-1}^j e^\tau h(\tau) d\tau + \int_n^{n+r} e^\tau h(\tau) d\tau \\
  &  \leq & \sum_{j=1}^n e^j \int_{j-1}^j h(\tau) d\tau + e^{n+r}\int_n^{n+r} h(\tau) d\tau \\
  & \leq & \sum_{j=1}^n e^j \int_{j-1}^j h(\tau) d\tau + e^{n+1}\int_n^{n+1} h(\tau) d\tau.
 \end{eqnarray*}
 Inequality (\ref{eq57}) implies
 \begin{eqnarray*}
  \int_{j-1}^j h(\tau) d\tau \leq a_1 + a_2 & , & j = 1, \dots , n+1 .
 \end{eqnarray*}
Therefore,
\begin{eqnarray*}
  e^{-t} \int_0^t e^\tau h(\tau) d\tau & = &  e^{-n-r}\int_0^{n+r} e^\tau h(\tau) d\tau =
  e^{-n-r}\sum_{j=1}^{n+1} e^j \int_{j-1}^j h(\tau) d\tau \\
  & \leq & (a_1+a_2) e^{-n-r} \sum_{j=1}^n e^{j} = (a_1+a_2) e^{-n-r} \frac{e^{n+2}-e}{e-1} \\
  & = & (a_1 + a_2) \frac{e^{2-r} - e^{-n-r+1}}{e-1} \leq (a_1 + a_2) \frac{e^{2}}{e-1}.
 \end{eqnarray*}
Thus,
\[
  \displaystyle \sup_{t \geq 0} e^{-t} \int_{0}^{t} e^{\tau} h (\tau ) d \tau
  \leq (a_1 + a_2) \frac{e^{2}}{e-1}  < \infty ,
  \]
  which proves the lemma.

\bibliography{error_estimates}
\bibliographystyle{amsplain}
\end{document}